\documentclass[12pt]{article}
\usepackage{e-jc}
\usepackage{bm,subeqnarray,eqsection,indent,cite}
\usepackage{latexsym}
\usepackage{amsmath}           
\usepackage{amsfonts}          
\usepackage{amssymb}           
\usepackage{url}

\begin{document}

\title{\vspace*{-2.5cm} Noncommutative determinants, \\
           Cauchy--Binet formulae, \\
           and Capelli-type identities  \\[5mm]
 \large\bf I.~Generalizations of the Capelli and Turnbull identities \\[0mm]}

\author{
  {\small Sergio Caracciolo}                       \\[-2mm]
  {\small\it Dipartimento di Fisica and INFN}      \\[-2mm]
  {\small\it Universit\`a degli Studi di Milano}   \\[-2mm]
  {\small\it via Celoria 16}                       \\[-2mm]
  {\small\it I-20133 Milano, ITALY}                \\[-2mm]
  {\small\tt Sergio.Caracciolo@mi.infn.it}         \\[-2mm]
  {\protect\makebox[5in]{\quad}}
   \\[-2mm]
  {\small Alan D.~Sokal\thanks{Also at Department of Mathematics,
           University College London, London WC1E 6BT, England.}}  \\[-2mm]
  {\small\it Department of Physics}       \\[-2mm]
  {\small\it New York University}         \\[-2mm]
  {\small\it 4 Washington Place}          \\[-2mm]
  {\small\it New York, NY 10003 USA}      \\[-2mm]
  {\small\tt sokal@nyu.edu}           \\[-2mm]
  {\protect\makebox[5in]{\quad}}  
   \\[-2mm]
  {\small Andrea Sportiello}                       \\[-2mm]
  {\small\it Dipartimento di Fisica and INFN}      \\[-2mm]
  {\small\it Universit\`a degli Studi di Milano}   \\[-2mm]
  {\small\it via Celoria 16}                       \\[-2mm]
  {\small\it I-20133 Milano, ITALY}                \\[-2mm]
  {\small\tt Andrea.Sportiello@mi.infn.it}         \\[-3mm]
   \\[-2mm]
}
\date{\dateline{September 20, 2008}{August 3, 2009}\\
\small Mathematics Subject Classification: 
15A15 (Primary);
05A19, 05A30, 05E15, 13A50, 15A24, 15A33, 15A72, 17B35, 20G05 (Secondary).}
\maketitle
\thispagestyle{empty}   

\vspace{-7mm}

\begin{abstract}
We prove, by simple manipulation of commutators,
two noncommutative generalizations of the Cauchy--Binet formula
for the determinant of a product.
As special cases we obtain elementary proofs of
the Capelli identity from classical invariant theory
and of Turnbull's Capelli-type identities for
symmetric and antisymmetric matrices.
\end{abstract}

\bigskip
\noindent
{\bf Key Words:}  Determinant, noncommutative determinant,
row-determinant, column-determinant, Cauchy--Binet theorem,
permanent, noncommutative ring, Capelli identity, Turnbull identity,
Cayley identity,
classical invariant theory, representation theory, Weyl algebra,
right-quantum matrix, Cartier--Foata matrix, Manin matrix.

\bigskip
\noindent

\clearpage

\newtheorem{defin}{Definition}[section]
\newtheorem{definition}[defin]{Definition}
\newtheorem{prop}[defin]{Proposition}
\newtheorem{proposition}[defin]{Proposition}
\newtheorem{lem}[defin]{Lemma}
\newtheorem{lemma}[defin]{Lemma}
\newtheorem{guess}[defin]{Conjecture}
\newtheorem{ques}[defin]{Question}
\newtheorem{question}[defin]{Question}
\newtheorem{prob}[defin]{Problem}
\newtheorem{problem}[defin]{Problem}
\newtheorem{thm}[defin]{Theorem}
\newtheorem{theorem}[defin]{Theorem}
\newtheorem{cor}[defin]{Corollary}
\newtheorem{corollary}[defin]{Corollary}
\newtheorem{conj}[defin]{Conjecture}
\newtheorem{conjecture}[defin]{Conjecture}
\newtheorem{examp}[defin]{Example}
\newtheorem{example}[defin]{Example}
\newtheorem{claim}[defin]{Claim}

\renewcommand{\theenumi}{\alph{enumi}}
\renewcommand{\labelenumi}{(\theenumi)}
\def\prf{\par\noindent{\bf Proof.\enspace}\rm}
\def\rmk{\par\medskip\noindent{\bf Remark.\enspace}\rm}

\newcommand{\be}{\begin{equation}}
\newcommand{\ee}{\end{equation}}
\newcommand{\<}{\langle}
\renewcommand{\>}{\rangle}
\newcommand{\widebar}{\overline}
\def\reff#1{(\protect\ref{#1})}
\def\spose#1{\hbox to 0pt{#1\hss}}
\def\ltapprox{\mathrel{\spose{\lower 3pt\hbox{$\mathchar"218$}}
 \raise 2.0pt\hbox{$\mathchar"13C$}}}
\def\gtapprox{\mathrel{\spose{\lower 3pt\hbox{$\mathchar"218$}}
 \raise 2.0pt\hbox{$\mathchar"13E$}}}
\def\textprime{${}^\prime$}
\def\proof{\par\medskip\noindent{\sc Proof.\ }}
\newcommand{\qed}{\quad $\Box$ \medskip \medskip}
\def\firstproof{\par\medskip\noindent{\sc First proof.\ }}
\def\secondproof{\par\medskip\noindent{\sc Second proof.\ }}
\def\proofof#1{\bigskip\noindent{\sc Proof of #1.\ }}
\def\alternateproofof#1{\bigskip\noindent{\sc Alternate proof of #1.\ }}
\def\half{ {1 \over 2} }
\def\third{ {1 \over 3} }
\def\twothird{ {2 \over 3} }
\def\smfrac#1#2{{\textstyle{#1\over #2}}}
\def\smsmfrac#1#2{{\scriptstyle{#1\over #2}}}
\def\smhalf{ {\smfrac{1}{2}} }
\def\smsmhalf{ {\smsmfrac{1}{2}} }
\newcommand{\real}{\mathop{\rm Re}\nolimits}
\renewcommand{\Re}{\mathop{\rm Re}\nolimits}
\newcommand{\imag}{\mathop{\rm Im}\nolimits}
\renewcommand{\Im}{\mathop{\rm Im}\nolimits}
\newcommand{\sgn}{\mathop{\rm sgn}\nolimits}
\newcommand\supp{\mathop{\rm supp}\nolimits}
\newcommand{\diag}{\mathop{\rm diag}\nolimits}
\newcommand{\pf}{\mathop{\rm pf}\nolimits}
\newcommand{\hf}{\mathop{\rm hf}\nolimits}
\newcommand{\tr}{\mathop{\rm tr}\nolimits}
\newcommand{\per}{\mathop{\rm per}\nolimits}
\newcommand{\adj}{\mathop{\rm adj}\nolimits}
\newcommand{\Res}{\mathop{\rm Res}\nolimits}
\newcommand{\rowdet}{\mathop{\hbox{\rm row-det}}\nolimits}
\newcommand{\coldet}{\mathop{\hbox{\rm col-det}}\nolimits}
\newcommand{\Mdet}{\mathop{\hbox{Mdet}}\nolimits}
\def\hboxscript#1{ {\hbox{\scriptsize\em #1}} }
\def\hboxrm#1{ {\hbox{\scriptsize\rm #1}} }
\newcommand{\rowper}{\mathop{\hbox{\rm row-per}}\nolimits}
\newcommand{\colper}{\mathop{\hbox{\rm col-per}}\nolimits}
\newcommand{\restrict}{\upharpoonright}
\renewcommand{\emptyset}{\varnothing}

\def\Z{{\mathbb Z}}
\def\ZZ{{\mathbb Z}}
\def\R{{\mathbb R}}
\def\C{{\mathbb C}}
\def\CC{{\mathbb C}}
\def\N{{\mathbb N}}
\def\NN{{\mathbb N}}
\def\Q{{\mathbb Q}}

\newcommand{\gl}{{\mathfrak{gl}}}
\newcommand{\ooo}{{\mathfrak{o}}}
\newcommand{\spsp}{{\mathfrak{sp}}}

\newcommand{\T}{{\rm T}}
\newcommand{\AT}{{A}^{\T}}
\newcommand{\BT}{{B}^{\T}}
\newcommand{\XT}{{X}^{\T}}

\newcommand{\scra}{{\mathcal{A}}}
\newcommand{\scrb}{{\mathcal{B}}}
\newcommand{\scrc}{{\mathcal{C}}}
\newcommand{\scrd}{{\mathcal{D}}}
\newcommand{\scre}{{\mathcal{E}}}
\newcommand{\scrf}{{\mathcal{F}}}
\newcommand{\scrg}{{\mathcal{G}}}
\newcommand{\scrh}{{\mathcal{H}}}
\newcommand{\scri}{{\mathcal{I}}}
\newcommand{\scrj}{{\mathcal{J}}}
\newcommand{\scrk}{{\mathcal{K}}}
\newcommand{\scrl}{{\mathcal{L}}}
\newcommand{\scrm}{{\mathcal{M}}}
\newcommand{\scrn}{{\mathcal{N}}}
\newcommand{\scro}{{\mathcal{O}}}
\newcommand{\scrp}{{\mathcal{P}}}
\newcommand{\scrq}{{\mathcal{Q}}}
\newcommand{\scrr}{{\mathcal{R}}}
\newcommand{\scrs}{{\mathcal{S}}}
\newcommand{\scrt}{{\mathcal{T}}}
\newcommand{\scru}{{\mathcal{U}}}
\newcommand{\scrv}{{\mathcal{V}}}
\newcommand{\scrw}{{\mathcal{W}}}
\newcommand{\scrx}{{\mathcal{X}}}
\newcommand{\scry}{{\mathcal{Y}}}
\newcommand{\scrz}{{\mathcal{Z}}}

\newcommand{\bgamma}{{\boldsymbol{\gamma}}}
\newcommand{\bsigma}{{\boldsymbol{\sigma}}}
\newcommand{\balpha}{{\boldsymbol{\alpha}}}
\newcommand{\bbeta}{{\boldsymbol{\beta}}}
\renewcommand{\pmod}[1]{\;({\rm mod}\:#1)}
\def\psibar{{\bar{\psi}}}
\def\etabar{{\bar{\eta}}}
\def\chibar{{\bar{\chi}}}
\def\xibar{{\bar{\xi}}}
\def\lambdabar{{\bar{\lambda}}}
\def\mubar{{\bar{\mu}}}
\def\varphibar{{\bar{\varphi}}}
\def\phibar{{\bar{\phi}}}
\def\cz{\overline{z}}


\newenvironment{sarray}{
	  \textfont0=\scriptfont0
	  \scriptfont0=\scriptscriptfont0
	  \textfont1=\scriptfont1
	  \scriptfont1=\scriptscriptfont1
	  \textfont2=\scriptfont2
	  \scriptfont2=\scriptscriptfont2
	  \textfont3=\scriptfont3
	  \scriptfont3=\scriptscriptfont3
	\renewcommand{\arraystretch}{0.7}
	\begin{array}{l}}{\end{array}}

\newenvironment{scarray}{
	  \textfont0=\scriptfont0
	  \scriptfont0=\scriptscriptfont0
	  \textfont1=\scriptfont1
	  \scriptfont1=\scriptscriptfont1
	  \textfont2=\scriptfont2
	  \scriptfont2=\scriptscriptfont2
	  \textfont3=\scriptfont3
	  \scriptfont3=\scriptscriptfont3
	\renewcommand{\arraystretch}{0.7}
	\begin{array}{c}}{\end{array}}

\newcommand{\bydef}{:=}
\newcommand{\defby}{=:}
\renewcommand{\implies}{\Longrightarrow}
\renewcommand{\binom}[2]{\left(#1\atop#2\right)}
\newcommand{\bigpartial}{\displaystyle\partial}

%
\newcommand{\ef}[1]{\, #1}     

\newcommand{\Reof}[1]{\mathfrak{Re}(#1)}
\newcommand{\Imof}[1]{\mathfrak{Im}(#1)}
\newcommand{\eval}[1]{\left\langle {#1} \right\rangle}
\newcommand{\leval}[1]{\langle {#1} \rangle}
\newcommand{\reval}[1]{\overline{#1}}

\newcommand{\sspan}{\mathrm{span}} 
\newcommand{\kker}{\mathrm{ker}}
\newcommand{\rrank}{\mathrm{rank}}

\newcommand{\bigast}[1]{\underset{#1}{\textrm{{\huge $\ast$}}}}

\newcommand{\dx}[1] {\mathrm{d}{#1}}
\newcommand{\dede}[1]{\frac{\partial}{\partial #1}}
\newcommand{\deenne}[2]{\frac{\partial^#2}{\partial #1 ^#2}}
\newcommand{\vett}[1]{#1}

\newcommand{\tinyfrac}[2] {\genfrac{}{}{}{1}{#1}{#2} }
\newcommand{\Lfrac}[2] {\genfrac{}{}{}{0}{#1}{#2} }

\newtheorem{ansatz}{Ansatz}[section]
\newtheorem{theor}{Theorem}
\newtheorem{coroll}{Corollary}

\clearpage

\section{Introduction}

Let $R$ be a commutative ring, and let $A = (a_{ij})_{i,j=1}^n$
be an $n \times n$ matrix with elements in $R$.
Define as usual the determinant
\begin{equation}
   \det A  \;\bydef\;
   \sum_{\sigma \in \scrs_n} \sgn(\sigma) \prod_{i=1}^n  a_{i \sigma(i)}
   \;.
 \label{def.det}
\end{equation}
One of the first things one learns about the determinant
is the {\em multiplicative property}\/:
\begin{equation}
   \det(AB)  \;=\; (\det A)(\det B)  \;.
 \label{eq.detAB}
\end{equation}
More generally, if $A$ and $B$ are $m \times n$ matrices,
and $I$ and $J$ are subsets of $[n] \bydef \{1,2,\ldots,n\}$
of cardinality $|I|=|J|=r$,
then one has the {\em Cauchy--Binet formula}\/:
\begin{subeqnarray}
   \det \: (\AT B)_{IJ}
   & = &
   \sum_{\begin{scarray}
            L \subseteq [m] \\
            |L| = r
          \end{scarray}}
      (\det \: (A^\T)_{IL}) (\det B_{L J})   \\[2mm]
   & = &
   \sum_{\begin{scarray}
            L \subseteq [m] \\
            |L| = r
          \end{scarray}}
      (\det A_{L I}) (\det B_{L J})
 \label{eq.Cauchy-Binet}
\end{subeqnarray}
where $M_{IJ}$ denotes the submatrix of $M$ with rows $I$ and columns $J$
(kept in their original order).

If one wants to generalize these formulae to matrices with elements in a
{\em noncommutative}\/ ring $R$, the first problem one encounters
is that the definition \reff{def.det} is ambiguous without an ordering
prescription for the product.
Rather, one can define numerous alternative ``determinants'':
for instance, the {\em column-determinant}\/
\begin{equation}
   \coldet A  \;\bydef\;
   \sum_{\sigma \in \scrs_n} \sgn(\sigma) \, 
   a_{\sigma(1) 1} \, a_{\sigma(2) 2} \,\cdots\, a_{\sigma(n) n}
 \label{def.coldet}
\end{equation}
and the {\em row-determinant}\/
\begin{equation}
   \rowdet A  \;\bydef\;
   \sum_{\sigma \in \scrs_n} \sgn(\sigma) \,
   a_{1 \sigma(1)} \, a_{2 \sigma(2)} \,\cdots\, a_{n \sigma(n)}
 \label{def.rowdet}
 \;.
\end{equation}
(Note that $\coldet A = \rowdet \AT$.)
Of course, in the absence of commutativity these ``determinants''
need not have all the usual properties of the determinant.

Our goal here is to prove the analogues of
\reff{eq.detAB}/\reff{eq.Cauchy-Binet}
for a fairly simple noncommutative case:
namely, that in which the elements of $A$
are in a suitable sense ``almost commutative'' among themselves (see below)
and/or the same for $B$,
while the commutators $[x,y] \bydef xy-yx$
of elements of $A$ with those of $B$
have the simple structure $[a_{ij}, b_{kl}] = - \delta_{ik} h_{jl}$.\footnote{
   The minus sign is inserted for future convenience.
   We remark that this formula makes sense even if the ring $R$
   lacks an identity element, as $\delta_{ik} h_{jl}$ is simply a shorthand
   for $h_{jl}$ if $i=k$ and 0 otherwise.
}
More precisely, we shall need the following type of commutativity
among the elements of $A$ and/or $B$:

\begin{definition}
  \label{def1.1}
Let $M = (M_{ij})$ be a (not-necessarily-square) matrix
with elements in a (not-necessarily-commutative) ring $R$.
$\!$Then we say that $M\!$ is {\em column-pseudo-commutative}\/ in case
\be
   [M_{ij}, M_{kl}]  \;=\; [M_{il}, M_{kj}]  \qquad\hbox{\rm for all } i,j,k,l
 \label{def.colpc.1}
\ee
and
\be
   [M_{ij}, M_{il}]  \;=\;  0   \qquad\hbox{\rm for all } i,j,l
   \;.
 \label{def.colpc.2}
\ee
We say that $M$ is {\em row-pseudo-commutative}\/ in case
$M^\T$ is column-pseudo-commutative.
\end{definition}

In Sections~\ref{sec.prelim} and \ref{sec.capelli}
we will explain the motivation for this strange definition,
and show that it really is the natural type of commutativity
for formulae of Cauchy--Binet type.\footnote{
   Similar notions arose already two decades ago
   in Manin's work on quantum groups \cite{Manin_87,Manin_88,Manin_91}.
   For this reason, some authors \cite{Chervov_08}
   call a row-pseudo-commutative matrix a {\em Manin matrix}\/;
   others \cite{Konvalinka_07,Konvalinka_08a,Konvalinka_08b}
   call it a {\em right-quantum matrix}\/.
   See the historical remarks at the end of Section~\ref{sec.prelim}.
}
Suffice it to observe now that column-pseudo-commutativity is
a fairly weak condition:
for instance, it is weaker than assuming that
$[M_{ij}, M_{kl}] = 0$ whenever $j \neq l$.
In many applications (though not all, see Example~\ref{example.itoh} below)
we will actually have $[a_{ij},a_{kl}] = [b_{ij},b_{kl}] = 0$
for {\em all}\/ $i,j,k,l$.
Note also that \reff{def.colpc.1} implies \reff{def.colpc.2}
if the ring $R$ has the property that $2x=0$ implies $x=0$.

The main result of this paper is the following:


\begin{proposition}[noncommutative Cauchy--Binet]
  \label{prop.capelli_rectangular}
{\quad}Let $R$ be a (not-necessarily-{\break}commutative) ring,
and let $A$ and $B$ be $m \times n$ matrices with elements in $R$. 
Suppose that
\be
      [a_{ij}, b_{kl}]   \;=\;  - \delta_{ik} h_{jl}
   \label{eq.hyp.prop.capelli_rectangular}
\ee
where $(h_{jl})_{j,l=1}^n$ are elements of $R$.
Then, for any $I, J \subseteq [n]$ of cardinality $|I|=|J|=r$:
\begin{itemize}
\item[(a)] If $A$ is column-pseudo-commutative, then
\be
   \sum_{\begin{scarray}
            L \subseteq [m] \\
            |L| = r
         \end{scarray}}
   (\coldet \: (\AT)_{IL}) (\coldet B_{L J})
     \;=\;   \coldet[(\AT B)_{I J} + Q_{\rm col}] 
 \label{eq.prop.capelli_rectangular.a}
\ee
where
\be
     (Q_{\rm col})_{\alpha\beta}  \;=\; (r-\beta)\, h_{i_\alpha j_\beta}
  \label{def.Qcol}
\ee
for $1\leq \alpha, \beta \leq r$.
\item[(b)] If $B$ is column-pseudo-commutative, then
\be
   \sum_{\begin{scarray}
            L \subseteq [m] \\
            |L| = r
         \end{scarray}}
   (\rowdet \: (\AT)_{IL}) (\rowdet B_{L J})
      \;=\;   \rowdet[(\AT B)_{I J} + Q_{\rm row}] 
 \label{eq.prop.capelli_rectangular.b}
\ee
where
\be
     (Q_{\rm row})_{\alpha\beta}  \;=\;  (\alpha-1)\, h_{i_\alpha j_\beta}  
 \label{def.Qrow}
\ee
for $1\leq \alpha, \beta \leq r$.
\end{itemize}
In particular,
\begin{itemize}
\item[(c)]  If $[a_{ij}, a_{kl}] = 0$ and $[b_{ij}, b_{kl}] = 0$
whenever $j \neq l$, then
\begin{subeqnarray}
   \sum_{\begin{scarray}
            L \subseteq [m] \\
            |L| = r
         \end{scarray}}
   (\det \: (\AT)_{I L}) (\det B_{L J})
   & = &
   \coldet[(\AT B)_{I J} + Q_{\rm col}]   \\[-7mm]
   & = &
   \rowdet[(\AT B)_{I J} + Q_{\rm row}]
 \label{eq.prop1.1c}
\end{subeqnarray}
\end{itemize}
\end{proposition}

\noindent
These identities can be viewed as a kind of ``quantum analogue''
of \reff{eq.Cauchy-Binet}, with the matrices $Q_{\rm col}$ and $Q_{\rm row}$
supplying the ``quantum correction''.
It is for this reason that we have chosen the letter $h$
to designate the matrix arising in the commutator.

Please note that the hypotheses of Proposition~\ref{prop.capelli_rectangular}
presuppose that $1 \le r \le n$
(otherwise $I$ and $J$ would be nonexistent or empty).
But $r > m$ is explicitly allowed:
in this case the left-hand side of
\reff{eq.prop.capelli_rectangular.a}/\reff{eq.prop.capelli_rectangular.b}/%
\reff{eq.prop1.1c}
is manifestly zero (since the sum over $L$ is empty),
but Proposition~\ref{prop.capelli_rectangular}
makes the nontrivial statement that the noncommutative determinant
on the right-hand side is also zero.

Note also that the hypothesis in part (c)
--- what we shall call {\em column-commutativity}\/,
 see Section~\ref{sec.prelim} ---
is sufficient to make the determinants of $A$ and $B$ well-defined
without any ordering prescription.
We have therefore written $\det$ (rather than $\coldet$ or $\rowdet$)
for these determinants.


Replacing $A$ and $B$ by their transposes
and interchanging $m$ with $n$ in Proposition~\ref{prop.capelli_rectangular},
we get the following ``dual'' version
in which the commutator $-\delta_{ik} h_{jl}$
is replaced by $-h_{ik} \delta_{jl}$:

\addtocounter{defin}{-1}
\begin{proposition}
\hspace*{-3mm} ${}^{\bf\prime}$ \hspace{1mm}
Let $R$ be a (not-necessarily-commutative) ring,
and let $A$ and $B$ be $m \times n$ matrices with elements in $R$. 
Suppose that
\be
      [a_{ij}, b_{kl}]   \;=\;  - h_{ik} \delta_{jl}
\ee
where $(h_{ik})_{i,k=1}^m$ are elements of $R$.
Then, for any $I, J \subseteq [m]$ of cardinality $|I|=|J|=r$:
\begin{itemize}
\item[(a)] If $A$ is row-pseudo-commutative, then
\be
   \sum_{\begin{scarray}
            L \subseteq [n] \\
            |L| = r
         \end{scarray}}
   (\coldet A_{IL}) (\coldet \: (\BT)_{L J})
     \;=\;   \coldet[(A \BT)_{I J} + Q_{\rm col}] 
\ee
where $Q_{\rm col}$ is defined in \reff{def.Qcol}.
\item[(b)] If $B$ is row-pseudo-commutative, then
\be
   \sum_{\begin{scarray}
            L \subseteq [n] \\
            |L| = r
         \end{scarray}}
   (\rowdet A_{IL}) (\rowdet \: (\BT)_{L J})
      \;=\;   \rowdet[(A \BT)_{I J} + Q_{\rm row}] 
\ee
where $Q_{\rm row}$ is defined in \reff{def.Qrow}.
\end{itemize}
In particular,
\begin{itemize}
\item[(c)]  If $[a_{ij}, a_{kl}] = 0$ and $[b_{ij}, b_{kl}] = 0$
whenever $i \neq k$, then
\begin{subeqnarray}
   \sum_{\begin{scarray}
            L \subseteq [n] \\
            |L| = r
         \end{scarray}}
   (\det A_{IL}) (\det \: (\BT)_{L J})
   & = &
   \coldet[(A \BT)_{I J} + Q_{\rm col}]   \\[-7mm]
   & = &
   \rowdet[(A \BT)_{I J} + Q_{\rm row}]
 \label{eq.prop1.1primec}
\end{subeqnarray}
\end{itemize}
\end{proposition}

When the commutator has the special form
$[a_{ij}, b_{kl}] = -h \delta_{ik} \delta_{jl}$,
then {\em both}\/ Propositions~\ref{prop.capelli_rectangular}
and \ref{prop.capelli_rectangular}${}'$ apply,
and by summing \reff{eq.prop1.1c}/\reff{eq.prop1.1primec}
over $I=J$ of cardinality $r$, we obtain:

\begin{corollary}
  \label{cor.capelli_rectangular}
Let $R$ be a (not-necessarily-commutative) ring,
and let $A$ and $B$ be $m \times n$ matrices with elements in $R$. 
Suppose that
\begin{subeqnarray}
      [a_{ij}, a_{kl}]   & = &   0  \\
      {}[b_{ij}, b_{kl}]   & = &   0  \\
      {}[a_{ij}, b_{kl}]   & = &  - h \delta_{ik} \delta_{jl}
 \label{eq.cor.capelli_rectangular}
\end{subeqnarray}
where $h \in R$.
Then, for any positive integer $r$, we have
\begin{subeqnarray}
   \sum_{\begin{scarray}
            I \subseteq [m] \\
            |I| = r
         \end{scarray}}
   \sum_{\begin{scarray}
            L \subseteq [n] \\
            |L| = r
         \end{scarray}}
   (\det A_{IL}) (\det B_{IL})
    & = &
   \sum_{\begin{scarray}
            I \subseteq [n] \\
            |I| = r
         \end{scarray}}
   \coldet[(\AT B)_{I I} + Q_{\rm col}]    \nonumber \\[-6mm] \\[2mm]
   & = &
   \sum_{\begin{scarray}
            I \subseteq [n] \\
            |I| = r
         \end{scarray}}
   \rowdet[(\AT B)_{I I} + Q_{\rm row}]    \nonumber \\[-6mm] \\[2mm]
    & = &
   \sum_{\begin{scarray}
            I \subseteq [m] \\
            |I| = r
         \end{scarray}}
   \coldet[(A \BT)_{I I} + Q_{\rm col}]  \nonumber \\[-6mm] \\[2mm]
   & = &
   \sum_{\begin{scarray}
            I \subseteq [m] \\
            |I| = r
         \end{scarray}}
   \rowdet[(A \BT)_{I I} + Q_{\rm row}]    \nonumber \\[-10mm]
 \label{eq.capelli}
\end{subeqnarray}
where
\begin{subeqnarray}
   Q_{\rm col}  & = &  h \diag(r-1,r-2,\ldots,0)  \\
   Q_{\rm row}  & = &  h \diag(0,1,\ldots,r-1)
 \label{eq.capelli.defQ}
\end{subeqnarray}
\end{corollary}

The cognoscenti will of course recognize
Corollary~\ref{cor.capelli_rectangular}
as (an abstract version of)
the Capelli identity \cite{Capelli_1887,Capelli_1888,Capelli_1890}
of classical invariant theory.
In Capelli's identity, the ring $R$ is the Weyl algebra $A_{m \times n}(K)$
over some field $K$ of characteristic 0 (e.g.\ $\Q$, $\R$ or $\C$)
generated by an $m \times n$ collection $X = (x_{ij})$
of commuting indeterminates (``positions'')
and the corresponding collection $\partial = (\partial/\partial x_{ij})$
of differential operators (proportional to ``momenta'');
we then take $A=X$ and $B=\partial$,
so that \reff{eq.cor.capelli_rectangular} holds with $h=1$.

The Capelli identity has a beautiful interpretation
in the theory of group representations \cite{Howe_91}:
Let $K = \R$ or $\C$, and consider the space $K^{m \times n}$
of $m \times n$ matrices with elements in $K$,
parametrized by coordinates $X = (x_{ij})$.
The group $GL(m) \times GL(n)$ acts on $K^{m \times n}$ by
\be
   (M,N) X  \;=\;  M^\T X N
 \label{def.GLmGLn_action}
\ee
where $M \in GL(m)$, $N \in GL(n)$ and $X \in K^{m \times n}$.
Then the infinitesimal action associated to \reff{def.GLmGLn_action}
gives a faithful representation
of the Lie algebra $\gl(m) \oplus \gl(n)$ 
by vector fields on $K^{m \times n}$ with linear coefficients:
\begin{subeqnarray}
  \gl(m) \colon\quad
      L_{ij} &\bydef& \sum_{l=1}^n x_{il} \, \frac{\partial}{\partial x_{jl}}
             \;=\; (X \partial^\T)_{ij}
             \qquad \hbox{for } 1 \le i,j \le m    \\[2mm]
  \gl(n) \colon\quad
      R_{ij} &\bydef& \sum_{l=1}^m x_{li} \, \frac{\partial}{\partial x_{lj}}
             \;=\; (X^\T \partial)_{ij}
             \qquad \hbox{for } 1 \le i,j \le n
\end{subeqnarray}
These vector fields have the commutation relations
\begin{subeqnarray}
   [L_{ij}, L_{kl}]   & = &  \delta_{jk} L_{il}  \,-\, \delta_{il} L_{kj} \\
   {}[R_{ij}, R_{kl}]   & = &  \delta_{jk} R_{il}  \,-\, \delta_{il} R_{kj} \\
   {}[L_{ij}, R_{kl}]   & = &  0
\end{subeqnarray}
characteristic of $\gl(m) \oplus \gl(n)$.
Furthermore, the action $(L,R)$ extends uniquely to a homomorphism
from the universal enveloping algebra $\scru(\gl(m) \oplus \gl(n))$
into the Weyl algebra $A_{m \times n}(K)$
[which is isomorphic to the algebra $\scrp\scrd(K^{m \times n})$
of polynomial-coefficient differential operators on $K^{m \times n}$].
As explained in \cite[secs.~1 and 11.1]{Howe_91},
it can be shown abstractly that any element of the Weyl algebra
that commutes with both $L$ and $R$ must be the image via $L$
of some element of the center of $\scru(\gl(m))$,
and also the image via $R$
of some element of the center of $\scru(\gl(n))$.
The Capelli identity \reff{eq.capelli} with $A=X$ and $B=\partial$
gives an explicit formula for the generators $\Gamma_r$
[$1 \le r \le \min(m,n)$] of this subalgebra,
from which it is manifest from (\ref{eq.capelli}a~or~b)
that $\Gamma_r$ belongs to the image under $R$ of $\scru(\gl(n))$
and commutes with the image under $L$ of $\scru(\gl(m))$,
and from (\ref{eq.capelli}c~or~d) the reverse fact.
See \cite{Fulton_91,Howe_89,Howe_91,Itoh_04,Kraft_96,Umeda_98,Umeda_08,Weyl_46}
for further discussion of the role of the Capelli identity
in classical invariant theory and representation theory,
as well as for proofs of the identity.

Let us remark that Proposition~\ref{prop.capelli_rectangular}${}'$
also contains Itoh's \cite{Itoh_04} Capelli-type identity
for the generators of the left action of $\ooo(m)$
on $m \times n$ matrices (see Example~\ref{example.itoh} below).

Let us also mention one important (and well-known) application of the
Capelli identity:  namely, it provides a simple proof
of the ``Cayley'' identity\footnote{
   The identity \reff{eq.intro.1}
   is conventionally attributed to Arthur Cayley (1821--1895);
   the generalization to arbitrary minors
   [see \reff{eq.cor.genCayley.2.minors} below]
   is sometimes attributed to Alfredo Capelli (1855--1910).
   The trouble is, neither of these formulae
   occurs anywhere --- as far as we can tell ---
   in the {\em Collected Papers}\/ of Cayley \cite{Cayley_collected}.
   Nor are we able to find these formulae in any of the relevant works
   of Capelli
   \cite{Capelli_1882,Capelli_1887,Capelli_1888,Capelli_1890,Capelli_02}.
   The operator $\Omega = \det(\partial)$ was indeed introduced
   by Cayley on the second page of his famous 1846 paper on invariants
   \cite{Cayley_1846};  it became known as Cayley's $\Omega$-process
   and went on to play an important role in classical invariant theory
   (see e.g.\ \cite{Weyl_46,Schur_68,Fulton_91,Kraft_96,Olver_99,Dolgachev_03}).
   But we strongly doubt that Cayley ever knew \reff{eq.intro.1}.
   See \cite{CSS_cayley,Abdesselam-Crilly-Sokal}
   for further historical discussion.
}
for $n \times n$ matrices,
\be
   \det(\partial) \, (\det X)^s  \;=\;
   s(s+1) \cdots (s+n-1) \, (\det X)^{s-1}
   \;.
 \label{eq.intro.1}
\ee
To derive \reff{eq.intro.1},
one simply applies both sides of the Capelli identity \reff{eq.capelli}
to $(\det X)^s$:
the ``polarization operators'' $L_{ij} = (X \partial^\T)_{ij}$
and $R_{ij} = (X^\T \partial)_{ij}$
act in a very simple way on $\det X$,
thereby allowing $\coldet(X \partial^\T + Q_{\rm col}) \, (\det X)^s$
and $\coldet(X^\T \partial + Q_{\rm col}) \, (\det X)^s$
to be computed easily;
they both yield $\det X$ times
the right-hand side of \reff{eq.intro.1}.\footnote{
   See e.g.\ \cite[p.~53]{Umeda_98} or \cite[pp.~569--570]{Howe_91}
   for derivations of this type.
}
In fact, by a similar method
we can use Proposition~\ref{prop.capelli_rectangular}
to prove a generalized ``Cayley'' identity
that lives in the Weyl algebra
(rather than just the polynomial algebra)
and from which the standard ``Cayley'' identity
can be derived as an immediate corollary:
see Proposition~\ref{prop.genCayley}
and Corollaries~\ref{cor.genCayley.1} and \ref{cor.genCayley.2}
in the Appendix.
See also \cite{CSS_cayley}
for alternate combinatorial proofs of a variety of Cayley-type identities.

Since the Capelli identity is widely viewed as ``mysterious''
\cite[p.~324]{Atiyah_73}
but also as a ``powerful formal instrument'' \cite[p.~39]{Weyl_46}
and a ``relatively deep formal result'' \cite[p.~40]{Stein_82},
it is of interest to provide simpler proofs.
Moreover, since the statement \reff{eq.capelli}/\reff{eq.capelli.defQ}
of the Capelli identity is purely algebraic/combinatorial,
it is of interest to give a purely algebraic/combinatorial proof,
independent of the apparatus of representation theory.
Such a combinatorial proof was given a decade ago
by Foata and Zeilberger \cite{Foata_94} for the case $m=n=r$,
but their argument was unfortunately somewhat intricate,
based on the construction of a sign-reversing involution.
The principal goal of the present paper is to provide
an extremely short and elementary algebraic proof of
Proposition~\ref{prop.capelli_rectangular}
and hence of the Capelli identity,
based on simple manipulation of commutators.
We give this proof in Section~\ref{sec.capelli}.

In 1948 Turnbull \cite{Turnbull_48} proved a Capelli-type identity
for {\em symmetric}\/ matrices (see also \cite{Wallace_53}),
and Foata and Zeilberger \cite{Foata_94} gave a combinatorial proof
of this identity as well.
Once again we prove a generalization:


\begin{proposition}[noncommutative Cauchy--Binet, symmetric version]
  \label{prop.capelli_symmetric}
Let $R$ be a (not-necessarily-commutative) ring,
and let $A$ and $B$ be $n \times n$ matrices with elements in $R$.
Suppose that
\begin{eqnarray}
   {}[a_{ij}, b_{kl}]  & = & -  h \, (\delta_{ik} \delta_{jl} +
                                      \delta_{il} \delta_{jk} )
\end{eqnarray}
where $h$ is an element of $R$.
\begin{itemize}
\item[(a)] Suppose that $A$ is column-pseudo-commutative and symmetric;
and if $n=2$, suppose further that either
\begin{itemize}
   \item[(i)] the ring $R$ has the property that $2x=0$ implies $x=0$, or
   \item[(ii)]  $[a_{12},h] = 0$.
\end{itemize}
Then, for any $I, J \subseteq [n]$ of cardinality $|I|=|J|=r$, we have
\begin{subeqnarray}
   \sum_{\begin{scarray}
            L \subseteq [n] \\
            |L| = r
         \end{scarray}}
   (\coldet A_{LI}) (\coldet B_{L J})
   & = &
   \coldet[(\AT B)_{I J} + Q_{\rm col}]   \\[-7mm]
   & = &
   \coldet[(A B)_{I J} + Q_{\rm col}]
\end{subeqnarray}
where
\be
     (Q_{\rm col})_{\alpha\beta}  \;=\; (r-\beta)\, h \delta_{i_\alpha j_\beta}
\ee
for $1\leq \alpha, \beta \leq r$.
\item[(b)] Suppose that $B$ is column-pseudo-commutative and symmetric;
and if $n=2$, suppose further that either
\begin{itemize}
   \item[(i)] the ring $R$ has the property that $2x=0$ implies $x=0$, or
   \item[(ii)]  $[b_{12},h] = 0$.
\end{itemize}
Then, for any $I, J \subseteq [n]$ of cardinality $|I|=|J|=r$, we have
\be
   \sum_{\begin{scarray}
            L \subseteq [n] \\
            |L| = r
         \end{scarray}}
   (\rowdet A_{LI}) (\rowdet B_{L J})
      \;=\;   \rowdet[(\AT B)_{I J} + Q_{\rm row}] 
\ee
where
\be
     (Q_{\rm row})_{\alpha\beta}  \;=\;
         (\alpha-1)\, h \delta_{i_\alpha j_\beta}  
\ee
for $1\leq \alpha, \beta \leq r$.
\end{itemize}
\end{proposition}

Turnbull \cite{Turnbull_48} and Foata--Zeilberger \cite{Foata_94}
proved their identity for a specific choice of matrices
$A = X^{\rm sym}$ and $B = \partial^{\rm sym}$ in a Weyl algebra,
but it is easy to see that their proof depends only on the
commutation properties and symmetry properties of $A$ and $B$.
Proposition~\ref{prop.capelli_symmetric} therefore generalizes
their work in three principal ways:
they consider only the case $r=n$, while we prove a general identity
of Cauchy--Binet type\footnote{
   See also Howe and Umeda \cite[sec.~11.2]{Howe_91}
   for a formula valid for general $r$,
   but involving a {\em sum}\/ over minors
   analogous to \reff{eq.capelli}.
};
they assume that {\em both}\/ $A$ and $B$ are symmetric,
while we show that it suffices for {\em one}\/ of the two to be symmetric;
and they assume that {\em both}\/ $[a_{ij},a_{kl}]=0$ and $[b_{ij},b_{kl}]=0$,
while we show that only {\em one}\/ of these plays any role
and that it moreover can be weakened to column-pseudo-commutativity.\footnote{
   This last weakening is, however, much less substantial
   than it might appear at first glance,
   because a matrix $M$ that is column-pseudo-commutative
   {\em and}\/ symmetric necessarily satisfies
   $2 [M_{ij},M_{kl}] = 0$ for {\em all}\/ $i,j,k,l$
   (see Lemma~\ref{lemma.commutators+symmetry} for the easy proof).
   In particular, in a ring $R$ in which $2x=0$ implies $x=0$,
   column-pseudo-commutativity plus symmetry implies full commutativity.
}
We prove Proposition~\ref{prop.capelli_symmetric}
in Section~\ref{sec.symmetric}.\footnote{
   In the first preprint version of this paper
   we mistakenly failed to include the extra hypotheses (i) or (ii)
   in Proposition~\ref{prop.capelli_symmetric} when $n=2$.
   For further discussion, see Section~\ref{sec.symmetric}
   and in particular Example~\ref{example.2by2.symmetric}.
}

Finally, Howe and Umeda \cite[eq.~(11.3.20)]{Howe_91}
and Kostant and Sahi \cite{Kostant_91}
independently discovered and proved
a Capelli-type identity for {\em antisymmetric}\/ matrices.\footnote{
   See also \cite{Kinoshita_02} for related work.
}
Unfortunately, Foata and Zeilberger \cite{Foata_94} were unable
to find a combinatorial proof of the Howe--Umeda--Kostant--Sahi identity;
and we too have been (thus far) unsuccessful.
We shall discuss this identity further in Section~\ref{sec.further}.

Both Turnbull \cite{Turnbull_48} and Foata--Zeilberger \cite{Foata_94}
also considered a different (and admittedly less interesting)
antisymmetric analogue of the Capelli identity,
which involves a generalization of the {\em permanent}\/ of a matrix $A$,
\be
   \per A  \;\bydef\;
   \sum_{\sigma \in \scrs_n} \prod_{i=1}^n  a_{i \sigma(i)}
   \;,
 \label{def.per}
\ee
to matrices with elements in a noncommutative ring $R$.
Since the definition \reff{def.per} is ambiguous
without an ordering prescription for the product,
we consider the {\em column-permanent}\/
\begin{equation}
   \colper A  \;\bydef\;
   \sum_{\sigma \in \scrs_n}  \, 
   a_{\sigma(1) 1} \, a_{\sigma(2) 2} \,\cdots\, a_{\sigma(n) n}
 \label{def.colper}
\end{equation}
and the {\em row-permanent}\/
\begin{equation}
   \rowper A  \;\bydef\;
   \sum_{\sigma \in \scrs_n}  \,
   a_{1 \sigma(1)} \, a_{2 \sigma(2)} \,\cdots\, a_{n \sigma(n)}
 \label{def.rowper}
 \;.
\end{equation}
(Note that $\colper A = \rowper \AT$.)
We then prove the following slight generalization
of Turnbull's formula:

\begin{proposition}[Turnbull's antisymmetric analogue]
  \label{prop.Turnbull-antisymmetric}
Let $R$ be a (not-necessarily-commutative) ring,
and let $A$ and $B$ be $n \times n$ matrices with elements in $R$.
Suppose that
\begin{eqnarray}
   {}[a_{ij}, b_{kl}]  & = & -  h \, (\delta_{ik} \delta_{jl} -
                                      \delta_{il} \delta_{jk} )
\end{eqnarray}
where $h$ is an element of $R$.
Then, for any $I, J \subseteq [n]$ of cardinality $|I|=|J|=r$:
\begin{itemize}
\item[(a)] If $A$ is antisymmetric off-diagonal
(i.e., $a_{ij} = -a_{ji}$ for $i \neq j$)
and $[a_{ij},h]=0$ for all $i,j$, we have
\begin{subeqnarray}
\lefteqn{ \hspace*{-3cm}  \sum_{\sigma\in{\cal S}_r} 
     \sum_{l_{1},\ldots, l_{r} \in [n]}
         a_{l_{1} i_{\sigma(1)}} \cdots a_{l_{r} i_{\sigma(r)}}
         b_{l_{1} j_{1}} \cdots b_{l_{r} j_{r}}\, =} \nonumber  \\
   & = &
   \colper[(\AT B)_{I J} - Q_{\rm col}]   \\[1mm]
   & = &
  (-1)^r\,  \colper[ (A B)_{I J} + Q_{\rm col}] 
 \label{eq.prop.Turnbull-antisymmetric.a}
\end{subeqnarray}
where
\be
     (Q_{\rm col})_{\alpha\beta}  \;=\; (r-\beta)\, h \delta_{i_\alpha j_\beta}
\ee
for $1\leq \alpha, \beta \leq r$.
\item[(b)] If $B$ is antisymmetric off-diagonal
(i.e., $b_{ij} = -b_{ji}$ for $i \neq j$)
and $[b_{ij},h]=0$ for all $i,j$, we have
\be
\sum_{\sigma\in{\cal S}_r} 
     \sum_{l_{1},\ldots, l_{r} \in [n]}
         a_{l_{1} i_{\sigma(1)}} \cdots a_{l_{r} i_{\sigma(r)}}
         b_{l_{1} j_{1}} \cdots b_{l_{r} j_{r}}
     \;=\;   \rowper[(\AT B)_{I J} - Q_{\rm row}] 
 \label{eq.prop.Turnbull-antisymmetric.b}
\ee
where
\be
     (Q_{\rm row})_{\alpha\beta}  \;=\;
         (\alpha-1)\, h \delta_{i_\alpha j_\beta}  
\ee
for $1\leq \alpha, \beta \leq r$.
\end{itemize}
\end{proposition}

\noindent
Note that no requirements are imposed on the $[a,a]$ and $[b,b]$ commutators
(but see the Remark at the end of Section~\ref{sec.symmetric}).

Let us remark that if $[a_{ij},b_{kl}]=0$, then the left-hand side of
\reff{eq.prop.Turnbull-antisymmetric.a}/\reff{eq.prop.Turnbull-antisymmetric.b}
is simply
\be
\sum_{\sigma\in{\cal S}_r} 
     \sum_{l_{1},\ldots, l_{r} \in [n]}
         a_{l_{1} i_{\sigma(1)}} \cdots a_{l_{r} i_{\sigma(r)}}
         b_{l_{1} j_{1}} \cdots b_{l_{r} j_{r}}
     \;=\;   \per (\AT B)_{I J}  \;,
\ee
so that Proposition~\ref{prop.Turnbull-antisymmetric}
becomes the {\em trivial}\/ statement
$\per (\AT B)_{I J} = \per (\AT B)_{I J}$.
So Turnbull's identity does {\em not}\/ reduce in the commutative case
to a formula of Cauchy--Binet type
--- indeed, no such formula exists for permanents\footnote{
   But see the Note Added at the end of this introduction.
}
---
which is why it is considerably less interesting than
the formulae of Cauchy--Binet--Capelli type for determinants.

Turnbull \cite{Turnbull_48} and Foata--Zeilberger \cite{Foata_94}
proved their identity for a specific choice of matrices
$A = X^{\rm antisym}$ and $B = \partial^{\rm antisym}$ in a Weyl algebra,
but their proof again depends only on the
commutation properties and symmetry properties of $A$ and $B$.
Proposition~\ref{prop.Turnbull-antisymmetric}
therefore generalizes their work in four principal ways:
they consider only the case $r=n$, while we prove a general identity
for minors;
they assume that {\em both}\/ $A$ and $B$ are antisymmetric,
while we show that it suffices for {\em one}\/ of the two to be antisymmetric
{\em plus an arbitrary diagonal matrix}\/;
and they assume that $[a_{ij},a_{kl}]=0$ and $[b_{ij},b_{kl}]=0$,
while we show that these commutators play no role.
We warn the reader that Foata--Zeilberger's \cite{Foata_94}
statement of this theorem contains a typographical error,
inserting a factor $\sgn(\sigma)$ that ought to be absent
(and hence inadvertently converting $\colper$ to $\coldet$).\footnote{
   Also, their verbal description of the other side of the identity ---
   ``the matrix product $X^\T P$ that appears on the right side
     of {\sc tur}${}'$ is taken with the assumption that the
     $x_{i,j}$ and $p_{i,j}$ commute'' --- is ambiguous,
   but we interpret it as meaning that all the factors $x_{i,j}$
   should be moved to the left, as is done
   on the left-hand side of
\reff{eq.prop.Turnbull-antisymmetric.a}/%
\reff{eq.prop.Turnbull-antisymmetric.b}.
}
We prove Proposition~\ref{prop.Turnbull-antisymmetric}
in Section~\ref{sec.symmetric}.\footnote{
   In the first preprint version of this paper
   we mistakenly failed to include the hypotheses that
   $[a_{ij},h] = 0$ or $[b_{ij},h]=0$.
   See the Remark at the end of Section~\ref{sec.symmetric}.
}

Finally, let us briefly mention some other generalizations of the
Capelli identity that have appeared in the literature.
One class of generalizations
\cite{Okounkov_96a,Okounkov_96b,Nazarov_98,Molev_07}
gives formulae for further elements in the (center of the)
universal enveloping algebra $\scru(\gl(n))$,
such as the so-called quantum immanants.
Another class of generalizations extends these formulae
to Lie algebras other than $\gl(n)$
\cite{Howe_91,Kostant_91,Kostant_93,Molev_99,Itoh_00,Itoh-Umeda_01,Itoh_04,%
Itoh_05,Wachi_06,Itoh_07,Molev_07}.
Finally, a third class of generalizations finds analogous formulae
in more general structures such as quantum groups \cite{Noumi_94,Noumi_96}
and Lie superalgebras \cite{Nazarov_97}.
Our approach is rather more elementary than all of these works:
we ignore the representation-theory context
and simply treat the Capelli identity as a noncommutative generalization of
the Cauchy--Binet formula.
A different generalization along vaguely similar lines
can be found in \cite{Mukhin_06}.


The plan of this paper is as follows:
In Section~\ref{sec.prelim} we make some preliminary comments
about the properties of column- and row-determinants.
In Section~\ref{sec.capelli} we prove
Propositions~\ref{prop.capelli_rectangular} and
\ref{prop.capelli_rectangular}${}'$
and Corollary~\ref{cor.capelli_rectangular}.
We also prove a variant of Proposition~\ref{prop.capelli_rectangular}
in which the hypothesis on the commutators $[a_{ij},a_{kl}]$
is weakened, at the price of a slightly weaker conclusion
(see Proposition~\ref{prop.capelli_rectangular_ANDREA}).
In Section~\ref{sec.symmetric} we prove
Propositions~\ref{prop.capelli_symmetric}
and~\ref{prop.Turnbull-antisymmetric}.
Finally, in Section~\ref{sec.further} we discuss whether
these results are susceptible of further generalization.
In the Appendix we prove a generalization of the ``Cayley'' identity
\reff{eq.intro.1}.

In a companion paper \cite{CS_capelli2}
we shall extend these identities to the (considerably more difficult) case
in which $[a_{ij}, b_{kl}] = - g_{ik} h_{jl}$
for general matrices $(g_{ik})$ and $(h_{jl})$,
whose elements do not necessarily commute.

\bigskip

{\em Note added.}\/  Subsequent to the posting of the present paper
in preprint form, Chervov, Falqui and Rubtsov \cite{Chervov_09}
posted an extremely interesting survey of the algebraic properties
of row-pseudo-commutative matrices (which they call ``Manin matrices'')
when the ring $R$ is an associative algebra over a field
of characteristic $\neq 2$.
In particular, Section~6 of \cite{Chervov_09} contains an interesting
generalization of the results of the present paper.\footnote{
   Chervov {\em et al.}\/ \cite{Chervov_09}
   also reformulated the hypotheses and proofs by using
   Grassmann variables (= exterior algebra) along the lines of
   \cite{Itoh-Umeda_01,Itoh_04}.
   This renders the proofs slightly more compact,
   and some readers may find that it renders the proofs
   more transparent as well (this is largely a question of taste).
   But we do think that the {\em hypotheses}\/ of the theorems
   are best stated without reference to Grassmann variables.
}
To state this generalization, note first that
the hypotheses of our Proposition~\ref{prop.capelli_rectangular}(a) are
\begin{itemize}
   \item[(i)] $A$ is column-pseudo-commutative, and
   \item[(ii)] $[a_{ij}, b_{kl}] \:=\: - \delta_{ik} h_{jl} \;.$
\end{itemize}
Left-multiplying (ii) by $a_{km}$ and summing over $k$, we obtain
\begin{itemize}
   \item[(ii${}'$)] $\sum\limits_k a_{km} \, [a_{ij}, b_{kl}]
                           \:+\: a_{im} h_{jl}
                     \:=\:
                     0  \;;$
\end{itemize}
moreover, the converse is true if $A$ is invertible.
Furthermore, (i) and (ii) imply
\begin{itemize}
   \item[(iii)]  $[a_{ij}, h_{ls}] \:=\: [a_{il}, h_{js}]$
\end{itemize}
as shown in Lemma~\ref{lemma.h_commute} below.
Then, Chervov {\em et al.}\/ \cite[Theorem~6]{Chervov_09} observed in essence
(translated back to our own language)
that our proof of Proposition~\ref{prop.capelli_rectangular}(a)
used only (i), (ii${}'$) and (iii),
and morover that (ii${}'$) can be weakened to\footnote{
   Here we have made the translations from their notation to ours
   ($M \to A^\T$, $Y \to B$, $Q \to H$)
   and written their hypotheses without reference to Grassmann variables.
   Their Conditions 1 and 2 then correspond to (ii${}''$) and (iii),
   respectively.
}
\begin{itemize}
   \item[(ii${}''$)] $\sum\limits_k a_{km} \, [a_{ij}, b_{kl}]
                           \:+\: a_{im} h_{jl}
                     \:=\:
                     [j \leftrightarrow m]$
\end{itemize}
--- that is, we need not demand the vanishing of
the left-hand side of (ii${}'$),
but merely of its antisymmetric part under $j \leftrightarrow m$,
{\em provided that we also assume (iii)}\/.
Their Theorem~6 also has the merit of including as a special case
not only Proposition~\ref{prop.capelli_rectangular}(a)
but also Proposition~\ref{prop.capelli_symmetric}.

Chervov {\em et al.}\/ \cite[Section~6.5]{Chervov_09} 
also provide an interesting rejoinder to our assertion above
that no formula of Cauchy--Binet type exists for permanents.
They show that if one defines a modified permanent
for submatrices involving {\em possibly repeated indices}\/,
which includes a factor $1/\nu!$ for each index that is repeated $\nu$ times,
then one obtains a formula of Cauchy--Binet type
in which the intermediate sum is over
$r$-tuples of {\em not necessarily distinct}\/ indices
$l_1 \le l_2 \le \ldots \le l_r$.
Moreover, this formula of Cauchy--Binet type
extends to a Capelli-type formula involving a ``quantum correction''
\cite[Theorems~11--13]{Chervov_09}.
In our opinion this is a very interesting observation,
which goes a long way to restore the analogy between
determinants and permanents
(and which in their formalism reflects the analogy between
 Grassmann algebra and the algebra of polynomials).


\section{Properties of column- and row-determinants}
  \label{sec.prelim}

In this section we shall make some preliminary observations
about the properties of column- and row-determinants,
stressing the following question:
Which commutation properties among the elements of the matrix
imply which of the standard properties of the determinant?
Readers who are impatient to get to the proof of our main results
can skim this section lightly.
We also call the reader's attention to the historical remarks
appended at the end of this section,
concerning the commutation hypotheses on matrix elements
that have been employed for theorems in noncommutative linear algebra.

Let us begin by recalling two elementary facts that we shall use repeatedly
in the proofs throughout this paper:

\begin{lemma}[Translation Lemma]
  \label{lemma.translation}
Let $\scra$ be an abelian group, and let $f \colon\, \scrs_n \to \scra$.
Then, for any $\tau \in \scrs_n$, we have
\be
   \sum_{\sigma \in \scrs_n}  \sgn(\sigma) \, f(\sigma)  \;=\;
   \sgn(\tau) \sum_{\sigma \in \scrs_n}  \sgn(\sigma) \, f(\sigma \circ \tau)
   \;.
\ee
\end{lemma}

\proof
Just note that both sides equal
$\sum\limits_{\sigma \in \scrs_n}
    \sgn(\sigma \circ \tau) \, f(\sigma \circ \tau)$.
\qed

\begin{lemma}[Involution Lemma]
  \label{lemma.involution}
Let $\scra$ be an abelian group, and let $f \colon\, \scrs_n \to \scra$.
Suppose that there exists a pair of distinct elements $i,j \in [n]$
such that
\be
     f(\sigma)  \;=\;  f(\sigma \circ (ij))
\ee
for all $\sigma \in \scrs_n$
[where $(ij)$ denotes the transposition interchanging $i$ with~$j$].
Then
\be
   \sum_{\sigma \in \scrs_n}  \sgn(\sigma) \, f(\sigma)  \;=\;   0
   \;.
\ee
\end{lemma}

\proof
We have
\begin{subeqnarray}
   \sum_{\sigma \in \scrs_n}  \sgn(\sigma) \, f(\sigma)
   & = &
   \sum_{\sigma \colon\, \sigma(i) < \sigma(j)}  \! \sgn(\sigma) \, f(\sigma)
   \;+\!\!
   \sum_{\sigma \colon\, \sigma(i) > \sigma(j)}  \! \sgn(\sigma) \, f(\sigma)
       \\[1mm] 
   & = &
   \sum_{\sigma \colon\, \sigma(i) < \sigma(j)}  \! \sgn(\sigma) \, f(\sigma)
   \;-\!\!
   \sum_{\sigma' \colon\, \sigma'(i) < \sigma'(j)}  \!\! \sgn(\sigma') \,
     f(\sigma' \circ (ij))
       \qquad \\[1mm] 
   & = &
   0  \;,
\end{subeqnarray}
where in the second line we made the change of variables
$\sigma' = \sigma \circ (ij)$ and used $\sgn(\sigma') = -\sgn(\sigma)$
[or equivalently used the Translation Lemma].
\qed

With these trivial preliminaries in hand,
let us consider noncommutative determinants.
Let $M = (M_{ij})$ be a matrix (not necessarily square)
with entries in a ring $R$.  Let us call $M$
\begin{itemize}
   \item  {\em commutative}\/ if $[M_{ij}, M_{kl}] = 0$ for all $i,j,k,l$;
   \item  {\em row-commutative}\/ if $[M_{ij}, M_{kl}] = 0$
                whenever $i \neq k$
                [i.e., all pairs of elements not in the same row commute];
   \item  {\em column-commutative}\/ if $[M_{ij}, M_{kl}] = 0$
                whenever $j \neq l$
                [i.e., all pairs of elements not in the same column commute];
   \item  {\em weakly commutative}\/ if $[M_{ij}, M_{kl}] = 0$
                whenever $i \neq k$ and $j \neq l$
                [i.e., all pairs of elements not in the same row or column
                 commute].
\end{itemize}
Clearly, if $M$ has one of these properties,
then so do all its submatrices $M_{IJ}$.
Also, $M$ is commutative if and only if
it is both row- and column-commutative.

Weak commutativity is a sufficient condition for the
determinant to be defined unambiguously without any ordering prescription,
since all the matrix elements in the product \reff{def.det}
differ in {\em both}\/ indices.
Furthermore, weak commutativity is sufficient for {\em single}\/
determinants to have most of their basic properties:

\begin{lemma}
   \label{lemma.weakQC.basic}
For weakly commutative square matrices:
\begin{itemize}
   \item[(a)] The determinant is antisymmetric under permutation of rows
     or columns.
   \item[(b)] The determinant of a matrix with two equal rows or columns
     is zero.
   \item[(c)] The determinant of a matrix equals the determinant of
     its transpose.
\end{itemize}
\end{lemma}

\noindent
The easy proof, which uses the Translation and Involution Lemmas,
is left to the reader (it is identical to the usual proof
in the commutative case).
We simply remark that
{\em if the ring $R$ has the property that $2x=0$ implies $x=0$}\/,
then antisymmetry under permutation of rows (or columns)
{\em implies}\/ the vanishing with equal rows (or columns).
But if the ring has elements $x \neq 0$ satisfying $2x=0$
(for instance, if $R = \Z_n$ for $n$ even),
then a slightly more careful argument,
using the Involution Lemma,
is needed to establish the vanishing with equal rows (or columns).

The situation changes, however, when we try to prove a formula
for the determinant of a product of {\em two}\/ matrices,
or more generally a formula of Cauchy--Binet type.
We are then inevitably led to consider products of matrix elements in which
some of the indices may be repeated ---
but only in {\em one}\/ of the two positions.
It therefore turns out (see Proposition~\ref{prop.easy_noncomm} below)
that we need something like {\em row}\/- or {\em column}\/-commutativity;
indeed, the result can be false without it
(see Example~\ref{example.easy_noncomm}).

Some analogues of Lemma~\ref{lemma.weakQC.basic}(a,b)
can nevertheless be obtained for the column- and row-determinants
under hypotheses {\em weaker}\/ than weak commutativity.
For brevity let us restrict attention to column-determinants;
the corresponding results for row-determinants
can be obtained by exchanging everywhere ``row'' with ``column''.

If $M = (M_{ij})_{i,j=1}^n$ is an $n \times n$ matrix
and $\tau \in \scrs_n$ is a permutation,
let us define the matrices obtained from $M$ by permutation of
rows or columns:
\begin{subeqnarray}
   ({}^\tau \! M)_{ij}  & \bydef &   M_{\tau(i) \, j}  \\[2mm]
   (M^\tau)_{ij}        & \bydef &   M_{i \, \tau(j)}
\end{subeqnarray}
We then have the following trivial result:

\begin{lemma}
   \label{lemma.coldet.rows}
For arbitrary square matrices:
\begin{itemize}
   \item[(a)] The column-determinant is antisymmetric under permutation of
     {\em rows}\/:
\be
   \coldet \, {}^\tau \! M
   \;=\;
   \sgn(\tau) \, \coldet M
\ee
for any permutation $\tau$.
   \item[(b)] The column-determinant of a matrix with two equal {\em rows}\/
     is zero.
\end{itemize}
\end{lemma}

\noindent
Indeed, statements (a) and (b) follow immediately from the
Translation Lemma and the Involution Lemma, respectively.


On the other hand, the column-determinant is {\em not}\/ in general
antisymmetric under permutation of {\em columns}\/,
nor is the column-determinant of a matrix with two equal columns
necessarily zero.  [For instance, in the Weyl algebra in one variable
over a field of characteristic $\neq 2$,
we have $\coldet \left( \!\! \displaystyle{ \begin{array}{cc}
                                                d & d \\
                                                x & x
                                             \end{array}
                                          }
                 \!\! \right) = dx-xd = 1$,
which is neither equal to $-1$ nor to 0.]
It is therefore natural to seek sufficient conditions for
these two properties to hold.  We now proceed to give a condition,
weaker than weak commutativity,
that entails the first property and {\em almost}\/ entails
the second property.

Let us begin by observing that
$\mu_{ijkl} \bydef [M_{ij}, M_{kl}]$
is manifestly antisymmetric under the simultaneous interchange
$i \leftrightarrow k$, $j \leftrightarrow l$.
So symmetry under one of these interchanges is equivalent to
antisymmetry under the other.  Let us therefore say that a matrix $M$ has
\begin{itemize}
   \item {\em row-symmetric}\/ (and {\em column-antisymmetric}\/)
      {\em commutators}\/
      if $[M_{ij}, M_{kl}] = [M_{kj}, M_{il}]$ for all $i,j,k,l$;
   \item {\em column-symmetric}\/ (and {\em row-antisymmetric}\/)
      {\em commutators}\/
      if $[M_{ij}, M_{kl}] = [M_{il}, M_{kj}]$ for all $i,j,k,l$.
\end{itemize}
Let us further introduce the same types of weakening
that we did for commutativity, saying that a matrix $M$ has
\begin{itemize}
  \item {\em weakly}\/ {\em row-symmetric}\/ (and {\em column-antisymmetric}\/)
      {\em commutators}\/
      if $[M_{ij}, M_{kl}] = [M_{kj}, M_{il}]$
      whenever $i \neq k$ and $j \neq l$;
  \item {\em weakly}\/ {\em column-symmetric}\/ (and {\em row-antisymmetric}\/)
      {\em commutators}\/
      if $[M_{ij}, M_{kl}] = [M_{il}, M_{kj}]$
      whenever $i \neq k$ and $j \neq l$.
\end{itemize}
(Note that row-symmetry is trivial when $i=k$,
 and column-symmetry is trivial when $j=l$.)
Obviously, each of these properties is inherited by
all the submatrices $M_{IJ}$ of $M$.
Also, each of these properties is manifestly weaker than the
corresponding type of commutativity.

The following fact is sometimes useful:

\begin{lemma}
   \label{lemma.commutators+symmetry}
Suppose that the square matrix $M$ has either row-symmetric or
column-symmetric commutators {\em and}\/
is either symmetric or antisymmetric.
Then $2 [M_{ij}, M_{kl}] = 0$ for all $i,j,k,l$.
In particular, if the ring $R$ has the property that $2x=0$ implies $x=0$,
then $M$ is commutative.
\end{lemma}

\proof
Suppose that $M$ has row-symmetric commutators
(the column-symmetric case is analogous)
and that $M^\T = \pm M$.
Then $[M_{ij}, M_{kl}] = [M_{kj}, M_{il}] = [M_{jk}, M_{li}]
      = [M_{lk}, M_{ji}] = [M_{kl}, M_{ij}]$,
where the first and third equalities use the row-symmetric commutators,
and the second and fourth equalities use symmetry or antisymmetry.
\qed

Returning to the properties of column-determinants, we have:

\begin{lemma}
 \label{lemma.coldet.columns}
If the square matrix $M$ has weakly row-symmetric commutators:
\begin{itemize}
   \item[(a)] The column-determinant is antisymmetric
under permutation of columns, i.e.
\be
   \coldet \, M^\tau
   \;=\;
   \sgn(\tau) \, \coldet M
\ee
for any permutation $\tau$.
   \item[(b)]  If $M$ has two equal columns, then $2 \coldet M = 0$.
(In particular, if $R$ is a ring in which $2x=0$ implies $x=0$,
 then $\coldet M = 0$.)
   \item[(c)]  If $M$ has two equal columns
   {\em and the elements in those columns commute among themselves}\/,
   then $\coldet M = 0$.
\end{itemize}
\end{lemma}

\proof
(a) It suffices to prove the claim when $\tau$ is the transposition
exchanging $i$ with $i+1$ (for arbitrary $i$).
We have
\begin{subeqnarray}
   \coldet M
   & = &
   \sum_{\sigma \in \scrs_n} \sgn(\sigma) \, 
      M_{\sigma(1),1} \, \cdots \,
      M_{\sigma(i),i} \, M_{\sigma(i+1),i+1} \,
      \, \cdots \, M_{\sigma(n),n}
   \\[2mm]
   & = &
   -\sum_{\sigma \in \scrs_n} \sgn(\sigma) \, 
      M_{\sigma(1),1} \, \cdots \, 
      M_{\sigma(i+1),i} \, M_{\sigma(i),i+1} \,
      \, \cdots \, M_{\sigma(n),n}
 \label{eq.coldetM}
\end{subeqnarray}
where the last equality uses the change of variables
$\sigma'= \sigma \circ (i,i+1)$ and the fact that
$\sgn(\sigma') = -\sgn(\sigma)$.
Similarly,
\begin{subeqnarray}
   \coldet M^\tau
   & = &
   \sum_{\sigma \in \scrs_n} \sgn(\sigma) \, 
      M_{\sigma(1),1} \, \cdots \, 
      M_{\sigma(i),i+1} \, M_{\sigma(i+1),i} \,
      \, \cdots \, M_{\sigma(n),n}
   \\[2mm]
   & = &
   -\sum_{\sigma \in \scrs_n} \sgn(\sigma) \, 
      M_{\sigma(1),1} \, \cdots \, 
      M_{\sigma(i+1),i+1} \, M_{\sigma(i),i} \,
      \, \cdots \, M_{\sigma(n),n}
   \;. \quad
 \label{eq.coldetMtau}
\end{subeqnarray}
It follows from (\ref{eq.coldetM}a) and (\ref{eq.coldetMtau}b) that
\begin{eqnarray}
 \coldet M \,+\, \coldet M^\tau
 & = &
 \sum_{\sigma \in \scrs_n} \sgn(\sigma) \,
      M_{\sigma(1),1} \, \cdots \,
      [M_{\sigma(i),i}, M_{\sigma(i+1),i+1}]
      \, \cdots \, M_{\sigma(n),n}
  \;.
  \nonumber \\[-3mm]
\end{eqnarray}
Under the hypothesis that $M$ has weakly row-symmetric commutators
[which applies here since $i \neq i+1$ and $\sigma(i) \neq \sigma(i+1)$],
the summand [excluding $\sgn(\sigma)$]
is invariant under $\sigma \mapsto \sigma \circ (i,i+1)$,
so the Involution Lemma implies that the sum is zero.

(b) is an immediate consequence of (a).

(c) Using (a), we may assume without loss of generality
that the two equal columns are adjacent (say, in positions 1 and 2).
Then, in
\be
   \coldet M  \;=\;
   \sum_{\sigma \in \scrs_n} \sgn(\sigma) \,
   M_{\sigma(1) 1} \, M_{\sigma(2) 2} \,\cdots\, M_{\sigma(n) n}
   \;,
 \label{eq.coldetM.c}
\ee
we have by hypothesis $M_{i1} = M_{i2}$ and
\be
   M_{\sigma(1) 1} \, M_{\sigma(2) 1}
   \;=\;
   M_{\sigma(2) 1} \, M_{\sigma(1) 1}
   \;,
\ee
so that the summand in \reff{eq.coldetM.c} [excluding $\sgn(\sigma)$]
is invariant under $\sigma \mapsto \sigma \circ (12)$;
the Involution Lemma then implies that the sum is zero.
\qed

The embarrassing factor of 2 in Lemma~\ref{lemma.coldet.columns}(b)
is not simply an artifact of the proof;
it is a fact of life when the ring $R$
has elements $x \neq 0$ satisfying $2x=0$:

\begin{example}
\rm
Let $R$ be the ring of $2 \times 2$ matrices with elements in
the field $GF(2)$, and let $\alpha$ and $\beta$ be any two noncommuting
elements of $R$
[for instance,
 $\alpha = \left( \!\! \displaystyle{ \begin{array}{cc}
                                          1 & 0 \\
                                          0 & 0
                                      \end{array}
                                    }
           \!\! \right)$
 and
 $\beta  = \left( \!\! \displaystyle{ \begin{array}{cc}
                                          0 & 1 \\
                                          1 & 0
                                      \end{array}
                                    }
           \!\! \right)$].
Then the matrix
$M =  \left( \!\! \displaystyle{ \begin{array}{cc}
                                     \alpha & \alpha \\
                                     \beta  & \beta
                                 \end{array}
                               }
      \!\! \right)$
has both row-symmetric and column-symmetric commutators
(and hence also row-{\em anti}\/symmetric
 and column-{\em anti}\/symmetric commutators! ---
 note that symmetry is {\em equivalent}\/ to antisymmetry
 in a ring of characteristic 2).
But $\coldet M = \alpha\beta-\beta\alpha \neq 0$.
\qed
\end{example}

In Proposition~\ref{prop.capelli_rectangular_ANDREA} below,
we shall prove a variant of Proposition~\ref{prop.capelli_rectangular}
that requires the matrix $\AT$ only to have row-symmetric commutators,
but at the price of multiplying everything by this embarrassing factor of 2.

If we want to avoid this factor of 2
by invoking Lemma~\ref{lemma.coldet.columns}(c),
then (as will be seen in Section~\ref{sec.capelli})
we shall need to impose a condition that is intermediate between
row-commutativity and row-symmetry:
namely, we say (as in Definition~\ref{def1.1}) that $M$ is
\begin{itemize}
   \item {\em row-pseudo-commutative}\/ if
      $[M_{ij},M_{kl}] = [M_{kj}, M_{il}]$ for all $i,j,k,l$ and
      $[M_{ij},M_{kj}] = 0$ for all $i,j,k$;
   \item {\em column-pseudo-commutative}\/ if
      $[M_{ij},M_{kl}] = [M_{il}, M_{jk}]$ for all $i,j,k,l$ and
      $[M_{ij},M_{il}] = 0$ for all $i,j,l$.
\end{itemize}
(Of course, the $[M,M]=[M,M]$ condition need be imposed
 only when $i \neq k$ and $j \neq l$,
 since in all other cases it is either trivial
 or else a consequence of the $[M,M]=0$ condition.)
We thus have $M$ row-commutative $\implies$
$M$ row-pseudo-commutative $\implies$
$M$ has row-symmetric commutators;
furthermore, the converse to the second implication holds
whenever $R$ is a ring in which $2x=0$ implies $x=0$.
Row-pseudo-commutativity thus turns out to be exactly
the strengthening of row-symmetry
that we need in order to apply Lemma~\ref{lemma.coldet.columns}(c)
and thus avoid the factor of 2
in Proposition~\ref{prop.capelli_rectangular_ANDREA},
i.e.\ to prove the full Proposition~\ref{prop.capelli_rectangular}.

The following intrinsic characterizations of
row-pseudo-commutativity and row-sym\-me\-try
are perhaps of some interest\footnote{
   Proposition~\ref{prop.intrinsic} is almost identical
   to a result of Chervov and Falqui \cite[Proposition~1]{Chervov_08},
   from whom we got the idea;
   but since they work in an associative algebra over a field of
   characteristic $\neq 2$, they don't need to distinguish
   between row-pseudo-commutativity and row-symmetry.
   They attribute this result to Manin \cite[top p.~199]{Manin_87}
   \cite{Manin_88,Manin_91}, but we are unable to find it there
   (or perhaps we have simply failed to understand what we have read).
   However, a result of similar flavor can be found in
   \cite[p.~193, Proposition]{Manin_87}
   \cite[pp.~7--8, Theorem 4]{Manin_88},
   and it is probably this to which the authors are referring.
}:

\begin{proposition}
  \label{prop.intrinsic}
Let $M = (M_{ij})$ be an $m \times n$ matrix with entries in a
(not-necessarily-commutative) ring $R$.
\begin{itemize}
   \item[(a)] Let $x_1, \ldots, x_n$ be commuting indeterminates,
and define for $1 \le i \le m$
the elements $\widetilde{x}_i = \sum\limits_{j=1}^n M_{ij} x_j$
in the polynomial ring $R[x_1,\ldots,x_n]$.
Then the matrix $M$ is row-pseudo-commutative
if and only if the elements $\widetilde{x}_1, \ldots, \widetilde{x}_m$
commute among themselves.
   \item[(b)] Let $\eta_1, \ldots, \eta_m$ be Grassmann indeterminates
(i.e.\ $\eta_i^2 = 0$ and $\eta_i \eta_j = -\eta_j \eta_i$),
and define for $1 \le j \le n$
the elements $\widetilde{\eta}_j = \sum\limits_{i=1}^m \!\eta_i M_{ij}$
in the Grassmann ring $R[\eta_1,\ldots,\eta_m]_{\rm Gr}$.
Then:
\begin{itemize}
   \item[(i)] The matrix $M$ has row-symmetric commutators
if and only if the elements $\widetilde{\eta}_1, \ldots, \widetilde{\eta}_n$
anticommute among themselves
(i.e.\ $\widetilde{\eta}_i \widetilde{\eta}_j =
        - \widetilde{\eta}_j \widetilde{\eta}_i$).
   \item[(ii)] The matrix $M$ is row-pseudo-commutative
if and only if the elements $\widetilde{\eta}_1, \ldots, \widetilde{\eta}_n$
satisfy all the Grassmann relations
$\widetilde{\eta}_i \widetilde{\eta}_j =
 - \widetilde{\eta}_j \widetilde{\eta}_i$
and $\widetilde{\eta}_i^2 = 0$.
\end{itemize}
\end{itemize}
\end{proposition}

\proof
(a) We have
\be
  [\widetilde{x}_i, \widetilde{x}_k]
  \;=\;
  \Biggl[ \sum\limits_j M_{ij} x_j ,\, \sum\limits_l M_{kl} x_l \Biggr]
  \;=\;
  \sum\limits_{j,l} [M_{ij}, M_{kl}] \, x_j x_l
  \;.
\ee
For $j \neq l$, the two terms in $x_j x_l = x_l x_j$ cancel
if and only if $[M_{ij}, M_{kl}] = -[M_{il}, M_{kj}]$;
and the latter equals $[M_{kj}, M_{il}]$.
For $j=l$, there is only one term, and it vanishes if and only if
$[M_{ij}, M_{kj}]=0$.

(b) We have
\be
   \widetilde{\eta}_j \widetilde{\eta}_l + \widetilde{\eta}_l \widetilde{\eta}_j
   \;=\;
   \sum\limits_{i,k} (\eta_i M_{ij} \eta_k M_{kl} + \eta_k M_{kl} \eta_i M_{ij})
   \;=\;
   \sum\limits_{i,k} \eta_i \eta_k \, [M_{ij}, M_{kl}]
\ee
since $\eta_k \eta_i = - \eta_i \eta_k$.
For $i \neq k$, the two terms in $\eta_i \eta_k = -\eta_k \eta_i$
cancel if and only if $[M_{ij}, M_{kl}] = [M_{kj}, M_{il}]$.
(Note that there is no term with $i=k$,
 so no further condition is imposed on the commutators $[M,M]$.)
On the other hand,
\be
   \widetilde{\eta}_j^2
   \;=\;
   \sum\limits_{i,k} \eta_i M_{ij} \eta_k M_{kj}
   \;=\;
   \sum\limits_{i < k} \eta_i \eta_k \, [M_{ij}, M_{kj}]
   \;,
\ee
which vanishes if and only if $[M_{ij}, M_{kj}] = 0$ for all $i,k$.
\qed

\bigskip

{\bf Some historical remarks.}
1.  Row-commutativity has arisen in some previous
work on noncommutative linear algebra,
beginning with the work of Cartier and Foata
on noncommutative extensions of the MacMahon master theorem
\cite[Th\'eor\`eme~5.1]{Cartier_69}.
For this reason, many authors
\cite{Konvalinka_07,Chervov_08,Konvalinka_08a,Konvalinka_08b}
call a row-commutative matrix a {\em Cartier--Foata matrix}\/.
See e.g.\ \cite{Cartier_69,Foata_79,Lalonde_96,Konvalinka_07,Konvalinka_08a}
for theorems of noncommutative linear algebra for row-commutative matrices;
and see also \cite[secs.~5 and 7]{Konvalinka_07}
for some beautiful $q$- and ${\bf q}$-generalizations.

2. Row-pseudo-commutativity has also arisen previously,
beginning (indirectly) with Manin's early work on quantum groups
\cite{Manin_87,Manin_88,Manin_91}.
Thus, some authors \cite{Chervov_08}
call a row-pseudo-commutative matrix a {\em Manin matrix}\/;
others \cite{Konvalinka_07,Konvalinka_08a,Konvalinka_08b}
call it a {\em right-quantum matrix}\/.
Results of noncommutative linear algebra
for row-pseudo-commutative matrices include
Cramer's rule for the inverse matrix \cite{Manin_88,Chervov_08,Konvalinka_08b}
and the Jacobi identity for cofactors \cite{Konvalinka_08b},
the formula for the determinant of block matrices \cite{Chervov_08},
Sylvester's determinantal identity \cite{Konvalinka_08a},
the Cauchy--Binet formula (Section~\ref{sec.capelli} below),
the Cayley--Hamilton theorem \cite{Chervov_08},
the Newton identities between $\tr M^k$ and coefficients of $\det(tI+M)$
\cite{Chervov_08},
and the MacMahon master theorem \cite{Konvalinka_07,Konvalinka_08b};
see also \cite[secs.~6 and 8]{Konvalinka_07} \cite{Konvalinka_08a,Konvalinka_08b}
for some beautiful $q$- and ${\bf q}$-generalizations.
See in particular \cite[Lemma~12.2]{Konvalinka_07}
for Lemma~\ref{lemma.coldet.columns}
specialized to row-pseudo-commutative matrices.

The aforementioned results suggest that
row-pseudo-commutativity is the natural hypothesis
for (most? all?)\ theorems of noncommutative linear algebra
involving the column-determinant.
Some of these results were derived earlier and/or have simpler proofs
under the stronger hypothesis of row-commutativity.

We thank Luigi Cantini for drawing our attention to the paper
\cite{Chervov_08}, from which we traced the other works cited here.

3.  Subsequent to the posting of the present paper
in preprint form, Chervov, Falqui and Rubtsov \cite{Chervov_09}
posted an extremely interesting survey of the algebraic properties
of row-pseudo-commutative matrices (which they call ``Manin matrices'')
when the ring $R$ is an associative algebra over a field
of characteristic $\neq 2$.
This survey discusses the results cited in \#2 above, plus many more;
in particular,
Section~6 of \cite{Chervov_09} contains an interesting
generalization of the results of the present paper
on Cauchy--Binet formulae and Capelli-type identities.
These authors state explicitly that
``the main aim of [their] paper is to argue the following claim:
linear algebra statements hold true for Manin matrices
in a form identical to the commutative case''
\cite[first sentence of Section~1.1]{Chervov_09}.

4.  The reader may well wonder (as one referee of the present paper did):
Since the literature already contains two competing terminologies
for the class of matrices in question
(``Manin'' and ``right-quantum''),
why muddy the waters by proposing yet another terminology
(``row-pseudo-commutative'') that is by no means guaranteed to catch on?
We would reply by stating our belief that a ``good'' terminology
ought to respect the symmetry $A \mapsto A^\T$;
or in other words, rows and columns ought to be treated on the same footing,
with neither one privileged over the other.
(For the same reason, we endeavor to treat
 the row-determinant and the column-determinant on an equal basis.)
We do not claim that our terminology is ideal ---
perhaps someone will find one that is more concise and easier to remember ---
but we do think that this symmetry property is important.

\section{Proof of the ordinary Capelli-type identities}
  \label{sec.capelli}

In this section we shall prove Proposition~\ref{prop.capelli_rectangular};
then Proposition~\ref{prop.capelli_rectangular}${}'$
and Corollary~\ref{cor.capelli_rectangular} follow immediately.
At the end we shall also prove a variant
(Proposition~\ref{prop.capelli_rectangular_ANDREA})
in which the hypotheses on the commutators are slightly weakened,
with a corresponding slight weakening of the conclusion.

It is convenient to begin by reviewing the proof of the
classical Cauchy--Binet formula \reff{eq.Cauchy-Binet}
where the ring $R$ is commutative.
First fix $L = \{l_1,\ldots,l_r\}$ with $l_1 < \ldots < l_r$,
and compute
\begin{subeqnarray}
   (\det \: (\AT)_{IL}) \, (\det B_{LJ})
   & = &
   \sum_{\tau,\pi \in{\cal S}_r} \sgn(\tau)\, \sgn(\pi)\,
        a_{l_{1} i_{\tau(1)}} \cdots a_{l_{r} i_{\tau(r)}}
        b_{l_{\pi(1)}j_{1}} \cdots b_{l_{\pi(r)} j_{r}}
   \nonumber \\[-4mm] \\
   &=&
   \sum_{\sigma,\pi \in{\cal S}_r} \sgn(\sigma)\,
        a_{l_{\pi(1)} i_{\sigma(1)}} \cdots
             a_{l_{\pi(r)} i_{\sigma(r)}}
        b_{l_{\pi(1)}j_{1}} \cdots b_{l_{\pi(r)} j_{r}}
   \;,
   \nonumber \\[-4mm]
 \label{eq.CBproof.1}
\end{subeqnarray}
where we have written $\sigma = \tau \circ \pi$
and exploited the commutativity of the elements of $A$ (but not of $B$).
Now the sum over $L$ and $\pi$ is equivalent to summing over
all $r$-tuples of {\em distinct}\/ elements $l_1,\ldots,l_r \in [m]$:
\be
   \sum_L (\det \: (\AT)_{IL}) \, (\det B_{LJ})
   \;=\;
   \!\!\!
   \sum_{l_1,\ldots,l_r \in [m] \,\hbox{\scriptsize\rm distinct}}
   \!\!\!
   f(l_1,\ldots,l_r)
   \: b_{l_{1} j_{1}} \cdots b_{l_{r} j_{r}}
   \;,
 \label{eq.sum_distinct}
\ee
where we have defined
\be
f(l_1,\ldots, l_r)  \;\bydef\;
    \sum_{\sigma\in{\cal S}_r} \sgn(\sigma)\,
    a_{l_{1} i_{\sigma(1)}} \cdots a_{l_{r} i_{\sigma(r)}}
 \label{def:f}
\ee
for arbitrary $l_1,\ldots,l_r \in [m]$.
Note now that $f(l_1,\ldots, l_r) = 0$
whenever two or more arguments take the same value,
because \reff{def:f} is then the determinant of a matrix
with two (or more) equal rows.
%
%
We can therefore add such terms to the sum \reff{eq.sum_distinct},
yielding
\begin{subeqnarray}
   \sum_L (\det \: (\AT)_{IL}) \, (\det B_{LJ})
   & = &
   \!\!\!
   \sum_{l_1,\ldots,l_r \in [m]}
       \sum_{\sigma\in{\cal S}_r} \sgn(\sigma)\,
       a_{l_{1} i_{\sigma(1)}} \cdots a_{l_{r} i_{\sigma(r)}}
       b_{l_{1} j_{1}} \cdots b_{l_{r} j_{r}}
   \nonumber \\[-4mm] \slabel{eq.CB.comm.a} \\
   & = &
   \sum_{\sigma\in{\cal S}_r} \sgn(\sigma) \,
        (\AT B)_{i_{\sigma(1)}j_{1}} \cdots (\AT B)_{i_{\sigma(r)}j_{r}}
   \\[1mm]
   & = &
   \det \: (\AT  B)_{IJ}
   \;,
 \label{eq.CBproof.2}
\end{subeqnarray}
where we have here commuted the $b$'s through the $a$'s.
Note that the order of the elements of $B$ remains unchanged
throughout these manipulations.

Let us also remark that this proof is valid even if $r > m$:
the starting sum \reff{eq.sum_distinct} is then empty,
since there do not exist
{\em distinct}\/ elements $l_1,\ldots,l_r \in [m]$;
but the sum \reff{eq.CB.comm.a} is nonempty,
since repetitions among the $l_1,\ldots,l_r$ are now allowed,
and we prove the nontrivial result that $\det \: (\AT  B)_{IJ} = 0$.
(Of course, in the commutative case this is no surprise,
since the matrix $\AT B$ has rank at most $m$;
but the corresponding noncommutative result will be less trivial.)

Now let us examine this proof more closely,
in order to see what commutation properties of the matrix elements
were really needed to make it work.
In the passage from (\ref{eq.CBproof.1}a) to (\ref{eq.CBproof.1}b),
the essence of the argument was that
\begin{subeqnarray}
   (\coldet \: (\AT)_{IL}) \, (\coldet B_{LJ})
   & = &
   \sum_{\pi \in{\cal S}_r} \sgn(\pi) \:
        (\coldet \: (\AT)_{IL}) \:
        b_{l_{\pi(1)}j_{1}} \cdots b_{l_{\pi(r)} j_{r}}
   \nonumber \\[-4mm] \\
   &=&
   \sum_{\pi \in{\cal S}_r} \sgn(\pi)^2 \:
        [\coldet \: ((\AT)_{IL})^\pi] \:
        b_{l_{\pi(1)}j_{1}} \cdots b_{l_{\pi(r)} j_{r}}
   \;,
   \nonumber \\[-4mm]
\end{subeqnarray}
where Lemma~\ref{lemma.coldet.columns}(a) justifies the passage
from the first line to the second;
so it suffices for $\AT$ to have weakly row-symmetric commutators.
In the argument that $f(l_1,\ldots, l_r) = 0$
whenever two or more arguments take the same value,
we need to apply Lemma~\ref{lemma.coldet.columns}(c)
to a matrix that is a submatrix of $\AT$
{\em with possibly repeated columns}\/;
therefore we need, in addition to weak row-symmetry,
the additional hypothesis that the matrix elements of $\AT$
within each column commute among themselves ---
or in other words, we need $\AT$ to be
row-pseudo-commutative (Definition~\ref{def1.1}).
Finally, in the step from (\ref{eq.CBproof.2}a) to (\ref{eq.CBproof.2}b),
we commuted the $b$'s through the $a$'s.
We have therefore proven:

\begin{proposition}[easy noncommutative Cauchy--Binet]
  \label{prop.easy_noncomm}
Let $R$ be a $\!$(not-nec\-es\-sar\-i\-ly-commutative) ring,
and let $A$ and $B$ be $m \times n$ matrices with elements in $R$.
Suppose that
\begin{itemize}
   \item[(a)]  $\AT$ is row-pseudo-commutative,
       i.e.\ $A$ is column-pseudo-commutative,
       i.e.\ $[a_{ij},a_{kl}] \!= [a_{il},a_{kj}]$
           whenever $i \neq k$ and $j \neq l$
       and $[a_{ij},a_{il}] = 0$ whenever $j \neq l$;
   \item[(b)]  the matrix elements of $A$ commute with those of $B$,
      i.e.\ $[a_{ij}, b_{kl}]  = 0$ for all $i,j,k,l$.
\end{itemize}
Then, for any $I, J \subseteq [n]$ of cardinality $|I|=|J|=r$, we have
\be
   \sum_{\begin{scarray}
            L \subseteq [m] \\
            |L| = r
         \end{scarray}}
   (\coldet \: (\AT)_{IL}) (\coldet B_{L J})
     \;=\;   \coldet \: (\AT B)_{I J}
   \;.
 \label{eq.prop.easy_noncomm}
\ee
\end{proposition}

\noindent
Note that no hypothesis whatsoever is needed concerning the commutators
$[b_{ij}, b_{kl}]$.

There is also a dual result using $\rowdet$,
in which $B$ is required to be column-pseudo-commutative
and no hypothesis is needed on the $[a,a]$ commutators.

\bigskip

The hypothesis in Proposition~\ref{prop.easy_noncomm}
that $A$ be column-pseudo-commutative really is necessary:

\begin{example}
   \label{example.easy_noncomm}
\rm
Let $\alpha$ and $\beta$ be any noncommuting elements of the ring $R$,
and let
$A = \left( \!\! \displaystyle{ \begin{array}{cc}
                                          \alpha & \beta \\
                                          0 & 0
                                      \end{array}
                                    }
           \!\! \right)$
and
$B  = \left( \!\! \displaystyle{ \begin{array}{cc}
                                          1 & 1 \\
                                          0 & 0
                                      \end{array}
                                    }
           \!\! \right)$
[let us assume for simplicity that the ring $R$ has an identity element],
so that
$\AT B = \left( \!\! \displaystyle{ \begin{array}{cc}
                                          \alpha & \alpha \\
                                          \beta  & \beta
                                      \end{array}
                                    }
           \!\! \right)$.
Then $A$ is row-commutative but not column-pseudo-commutative,
while the elements of $B$ commute with everything.
We have $\det \AT = \det B = 0$ but
$\coldet(\AT B) = \alpha\beta-\beta\alpha \neq 0$.

This example can be streamlined by dropping the second row
of the matrices $A$ and $B$,
i.e.\ considering it as an example with $m=1$, $n=2$ and $r=2$.
Then the left-hand side of \reff{eq.prop.easy_noncomm}
is an empty sum (since $r > m$),
but the right-hand side does not vanish.
\qed
\end{example}

\begin{example}
   \label{example.2by2}
\rm
It is instructive to consider the general case of $2 \times 2$ matrices
(i.e.\ $m=n=2$)
under the sole hypothesis that $[a_{ij}, b_{kl}] = 0$ for all $i,j,k,l$.
We have
\begin{eqnarray}
 \coldet(\AT B) \,-\, (\coldet \AT )(\coldet B)
  & = &
  \bigl( [a_{21},a_{12}] \,+\, [a_{11},a_{22}] \bigr) \, b_{21} b_{12}
       \nonumber \\[1mm]
  & &  \!\!\!\!\!\!+\;  [a_{11},a_{12}] \, b_{11} b_{12}  \;+\;
            [a_{21},a_{22}] \, b_{21} b_{22}
  \,, \qquad
 \label{eq.ab22}
\end{eqnarray}
where the terms on the first (resp.\ second) line of the right-hand side
come from the first (resp.\ second) step of the proof.
We see that column-pseudo-commutativity of $A$
is {\em precisely}\/ what we need in order to guarantee that
\reff{eq.ab22} vanishes for arbitrary matrices $B$.
\qed
\end{example}

\bigskip

We are now ready to consider Proposition~\ref{prop.capelli_rectangular},
which generalizes Proposition~\ref{prop.easy_noncomm}
by allowing nonzero commutators $[a_{ij}, b_{kl}] = - \delta_{ik} h_{jl}$,
thereby producing a ``quantum correction'' on the right-hand side
of the identity.
In the proof of Proposition~\ref{prop.capelli_rectangular}(a)
it will be necessary (as we shall see) to commute the $h$'s through the $a$'s.
We therefore begin with a lemma
giving an important property of such commutators:

\begin{lemma}
    \label{lemma.h_commute}
Let $R$ be a (not-necessarily-commutative) ring,
and let $A$ and $B$ be $m \times n$ matrices with elements in $R$.
Suppose that for all $i,k \in [m]$ and $j,l \in [n]$ we have
\begin{subeqnarray}
   [a_{ij}, a_{il}]   & = &    0      \\
   {}[a_{ij}, b_{kl}] & = &   - \delta_{ik} h_{jl}
 \label{eq.lemma.h_commute.1}
\end{subeqnarray}
where $(h_{jl})_{j,l=1}^n$ are elements of $R$.
Then, for all $i \in [m]$ and $j,l,s \in [n]$ we have
\be
  [a_{ij}, h_{ls}]
    \;=\;
  [a_{il}, h_{js}]
       \;.
 \label{eq.lemma.h_commute.2}
\ee
\end{lemma}

\noindent
Note the very weak hypothesis here on the $[a,a]$ commutators:
we require $[a_{ij}, a_{kl}] = 0$ only when $i=k$,
i.e.\ between different columns {\em within the same row}\/.
This is much weaker than the column-pseudo-commutativity
assumed in Proposition~\ref{prop.capelli_rectangular}(a),
as it imposes \reff{def.colpc.2} but omits \reff{def.colpc.1}.

\proof
For any indices $i,k,r \in [m]$ and $j,l,s \in [n]$,
we have the Jacobi identity
\be
    [a_{ij}, [a_{kl}, b_{rs}]] \,+\,
    [a_{kl}, [b_{rs}, a_{ij}]] \,+\,
    [b_{rs}, [a_{ij}, a_{kl}]]   \;=\;   0
    \;.
 \label{Jacobi}
\ee
By taking $k=r=i$ and using the hypotheses \reff{eq.lemma.h_commute.1},
we obtain the conclusion \reff{eq.lemma.h_commute.2}.
\qed

\medskip
\noindent
{\bf Remark.}
Since $A$ and $B$ play symmetrical roles in this problem
(modulo the substitution $h \mapsto -h^\T$),
a similar argument shows that if $[b_{ij}, b_{il}] = 0$,
then $[b_{ij}, h_{ls}] = [b_{is}, h_{lj}]$.
This will be relevant for Proposition~\ref{prop.capelli_rectangular}(b).
\qed

\bigskip

One consequence of Lemma~\ref{lemma.h_commute}
is that $h$ can be commuted through $a$ when it arises
inside a sum over permutations with the factor $\sgn(\sigma)$:

\begin{corollary}
  \label{cor.h_commute}
Fix distinct elements $\alpha,\beta \in [r]$
and fix a set $I \subseteq [n]$ of cardinality $|I|=r$.
Then, under the hypotheses of Lemma~\ref{lemma.h_commute}, we have
\be
   \sum\limits_{\sigma \in \scrs_r}
      \sgn(\sigma) \:
      F\!\left( \{\sigma(j)\}_{j \neq \alpha,\beta}^{\vphantom{X}} \right)
      [a_{l i_{\sigma(\alpha)}}, h_{i_{\sigma(\beta)} k}] \:
      G\!\left( \{\sigma(j)\}_{j \neq \alpha,\beta}^{\vphantom{X}} \right)
   \;=\;  0
 \label{eq.cor.h_commute}
\ee
for arbitrary functions $F,G \colon\, [r]^{r-2} \to R$
and arbitrary indices $l \in [m]$ and $k \in [n]$.
\end{corollary}

\proof
By Lemma~\ref{lemma.h_commute} we have
\be
   [a_{l i_{\sigma(\alpha)}}, h_{i_{\sigma(\beta)} j}]
   \;=\;
   [a_{l i_{\sigma(\beta)}}, h_{i_{\sigma(\alpha)} j}]
   \;.
\ee
This means that the summand in \reff{eq.cor.h_commute}
[excluding the factor $\sgn(\sigma)$]
is invariant under $\sigma \mapsto \sigma \circ (\alpha\beta)$.
The claim then follows immediately from the Involution Lemma.
\qed

We also have a dual version of Corollary~\ref{cor.h_commute},
along the lines of the Remark above,
stating that if $[b_{ij}, b_{il}] = 0$,
then sums involving $[b_{l j_{\sigma(\alpha)}}, h_{k j_{\sigma(\beta)}}]$
vanish.
Let us call this Corollary~\ref{cor.h_commute}${}'$.

\bigskip

We are now ready to prove Proposition~\ref{prop.capelli_rectangular}:

\proofof{Proposition~\ref{prop.capelli_rectangular}}
We begin with part (a).
The first two steps in the proof are identical to those in
Proposition~\ref{prop.easy_noncomm}:  we therefore have
\be
 \sum_{L} (\coldet \: (\AT)_{IL}) \, (\coldet B_{LJ}) 
 \;=\;
 \sum_{\sigma\in{\cal S}_r} \sgn(\sigma) \sum_{l_{1},\cdots, l_{r} \in [m]}
      a_{l_{1} i_{\sigma(1)}} \cdots a_{l_{r} i_{\sigma(r)}}
      b_{l_{1} j_{1}} \cdots b_{l_{r} j_{r}} 
 \;.
 \label{eq.mainproof.1}
\ee
It is only now that we have to work harder,
because of the noncommutativity of the $b$'s with the $a$'s.
Let us begin by moving the factor $b_{l_{1} j_{1}}$ to the left
until it lies just to the right of $a_{l_{1} i_{\sigma(1)}}$,
using the general formula
\be
   x_1 \, [x_2 \cdots x_r, y]
   \;=\;
   x_1 \sum_{s=2}^r  x_2 \cdots x_{s-1} \, [x_s,y] \, x_{s+1} \cdots x_r
\ee
with $x_\alpha = a_{l_\alpha i_{\sigma(\alpha)}}$ and $y = b_{l_{1} j_{1}}$.
This gives
\begin{eqnarray}
\lefteqn{  \sum_{\sigma\in{\cal S}_r} \sgn(\sigma)
     \sum_{l_{1},\ldots, l_{r} \in [m]}
         a_{l_{1} i_{\sigma(1)}} \cdots a_{l_{r} i_{\sigma(r)}}
         b_{l_{1} j_{1}} \cdots b_{l_{r} j_{r}}
        }   \nonumber \\[2mm]
& = &
   \sum_{\sigma\in{\cal S}_r} \sgn(\sigma)
   \sum_{l_{1},\ldots, l_{r} \in [m]}
         a_{l_{1}i_{\sigma(1)}}
         \left[  \vphantom{ \sum_{k=2}^{m}} 
                  b_{l_{1}j_{1}}
                  a_{l_{2}i_{\sigma(2)}}  \cdots a_{l_{r} i_{\sigma(r)}}
         \right.  \nonumber \\[-1mm]
&  & \,\left.   - \sum_{s=2}^{r}  \delta_{l_{1}l_{s}} \,
            a_{l_{2}i_{\sigma(2)}} \cdots a_{l_{s-1}i_{\sigma(s-1)}}
            h_{i_{\sigma(s)}j_{1}}
            a_{l_{s+1}i_{\sigma(s+1)}} \cdots a_{l_{r}i_{\sigma(r)}}
       \right]  b_{l_{2} j_{2}} \cdots b_{l_{r} j_{r}}
   \,.   \qquad
\end{eqnarray}
Now we repeatedly use Corollary~\ref{cor.h_commute}
to push the factor $h_{i_{\sigma(s)}j_{1}}$ to the left:  we obtain
\begin{subeqnarray}
\lefteqn{   \sum_{\sigma\in{\cal S}_r} \sgn(\sigma)
   \sum_{l_{1},\ldots, l_{r} \in [m]}
         a_{l_{1}i_{\sigma(1)}}
         \left[  \vphantom{ \sum_{k=2}^{m}} 
                  b_{l_{1}j_{1}}
                  a_{l_{2}i_{\sigma(2)}}  \cdots a_{l_{r} i_{\sigma(r)}}
         \right.
        }    \nonumber \\[-1mm]
&  & \,\left.   - \, \sum_{s=2}^{r} h_{i_{\sigma(s)}j_{1}} \delta_{l_{1}l_{s}}
            a_{l_{2}i_{\sigma(2)}} \cdots a_{l_{s-1}i_{\sigma(s-1)}}
            a_{l_{s+1}i_{\sigma(s+1)}} \cdots a_{l_{r}i_{\sigma(r)}}
       \right]  b_{l_{2} j_{2}} \cdots b_{l_{r} j_{r}}   \quad \\[2mm]
&=&
    \sum_{\sigma\in{\cal S}_r} \sgn(\sigma)
    \sum_{l_{2},\ldots, l_{r} \in [m]}
          \left[ \vphantom{ \sum_{k=2}^{n}}
                 (\AT B)_{i_{\sigma(1)} j_{1}}
                 a_{l_{2}i_{\sigma(2)}}  \cdots a_{l_{r}i_{\sigma(r)}}
          \right.   \nonumber \\[-1mm]
&&  \,\left.  - \,\sum_{s=2}^{r} h_{i_{\sigma(s)} j_{1}}
            a_{l_{2}i_{\sigma(2)}} \cdots a_{l_{s-1}i_{\sigma(s-1)}}
            a_{l_{s}i_{\sigma(1)}}
            a_{l_{s+1}i_{\sigma(s+1)}} \cdots a_{l_{r}i_{\sigma(r)}}
      \!\right]
            b_{l_{2} j_{2}} \cdots b_{l_{r} j_{r}} \quad \\[2mm]
&=&
    \sum_{\sigma\in{\cal S}_r} \sgn(\sigma)
    \left[ (\AT B)_{i_{\sigma(1)}j_{1}}  +
           \,\sum_{s=2}^{r} h_{i_{\sigma(1)}j_{1}} \right]  
    \sum_{l_{2},\ldots, l_{r} \in [m]}
            a_{l_{2}i_{\sigma(2)}} \cdots a_{l_{r}i_{\sigma(r)}}
            b_{l_{2} j_{2}} \cdots b_{l_{r} j_{r}}  \nonumber \\[-3mm] \\[2mm]
&=&
    \sum_{\sigma\in{\cal S}_r} \sgn(\sigma)
    \left[\AT B  + \,(r-1)\,h \right]_{i_{\sigma(1)}j_{1}}
    \!\!\sum_{l_{2},\ldots, l_{r} \in [m]}
         \!\!\!  a_{l_{2}i_{\sigma(2)}} \cdots a_{l_{r}i_{\sigma(r)}}
            b_{l_{2} j_{2}} \cdots b_{l_{r} j_{r}}
    \,, \qquad
\end{subeqnarray}
where we have simply executed the sum over $l_{1}$ and,
in the second summand, interchanged $\sigma(1)$ with $\sigma(s)$
[which multiplies $\sgn(\sigma)$ by $-1$].
This procedure can be now iterated to obtain
\begin{subeqnarray}
 \lefteqn{ \sum_{L} (\det \: (\AT)_{IL}) (\coldet B_{LJ})   }  \\[1mm]
 & = &
 \sum_{\sigma\in{\cal S}_r} \sgn(\sigma)
    \left[\AT B  + \,(r-1)\,h \right]_{i_{\sigma(1)}j_{1}}
    \left[\AT B  + \,(r-2)\,h \right]_{i_{\sigma(2)}j_{2}}
    \cdots
    \left[\AT B \right]_{i_{\sigma(r)}j_{r}}     \nonumber\\[-3mm] \\[1mm]
 & = &
 \coldet  \left[(\AT B)_{I J} + Q_{\rm col} \right]
   \,,
 \label{eq.mainproof.5}
\end{subeqnarray}
which is the desired result of part (a).

For part (b), let us start as before:
\begin{subeqnarray}
\lefteqn{ \sum_{L} (\rowdet (\AT)_{IL}) \, \left(\rowdet B_{LJ}\right) } \\
&=& \sum_{L}\sum_{\tau,\pi\in{\cal S}_r} \sgn(\tau) \sgn(\pi)\, 
        a_{l_{\tau(1)}i_{1}} \cdots a_{l_{\tau(r)} i_{r}}
        b_{l_{1} j_{\pi(1)}} \cdots b_{l_{r} j_{\pi(r)}}  \\
&=& \sum_{L}\sum_{\tau,\sigma\in{\cal S}_r} \sgn(\sigma)\,
        a_{l_{\tau(1)}i_{1}} \cdots a_{l_{\tau(r)} i_{r}}
        b_{l_{\tau(1)}j_{\sigma(1)}} \cdots b_{l_{\tau(r)} j_{\sigma(r)}}
\end{subeqnarray}
where we have written $\sigma = \pi \circ \tau$
and exploited the commutativity of the elements of $B$ (but not of $A$).
An argument as in Proposition~\ref{prop.easy_noncomm}
allows us to rewrite this as
\be
   \sum_{\sigma\in{\cal S}_r} \sgn(\sigma)
   \sum_{l_{1},\cdots, l_{r} \in [m]}
      a_{l_{1}i_{1}} \cdots a_{l_{r} i_{r}}
      b_{l_{1}j_{\sigma(1)}} \cdots b_{l_{r} j_{\sigma(r)}} 
   \;.
\ee
We first move the factor $a_{l_{r} i_{r}}$ to the right, giving
\begin{eqnarray}
& & \sum_{\sigma\in{\cal S}_r} \sgn(\sigma)
      \sum_{l_{1},\ldots, l_{r} \in [m]}
          a_{l_{1}i_{1}} \cdots a_{l_{r-1} i_{r-1}}
          \left[ \vphantom{\sum_{s=1}^{r-1}}
                 b_{l_{1}j_{\sigma(1)}} \cdots b_{l_{r-1} j_{\sigma(r-1)}}
                 a_{l_{r} i_{r}}  \right.    \nonumber \\[-1mm]
& & \left. \quad - \sum_{s=1}^{r-1}   \delta_{l_r l_s}
                b_{l_{1} j_{\sigma(1)}}  \cdots b_{l_{s-1} j_{\sigma(s-1)}}
                h_{i_r j_{\sigma(s)}}
                b_{l_{s+1} j_{\sigma(s+1)}} \cdots b_{l_{r-1} j_{\sigma(r-1)}}
    \right]   b_{l_{r} j_{\sigma(r)}}
    \;. \qquad
\end{eqnarray}
Now we repeatedly use Corollary~\ref{cor.h_commute}${}'$
to push the factor $h_{i_r j_{\sigma(s)}}$ to the right:  we obtain
\begin{subeqnarray}
&  & \sum_{\sigma\in{\cal S}_r} \sgn(\sigma)
      \sum_{l_{1},\ldots, l_{r} \in [m]}
          a_{l_{1}i_{1}} \cdots a_{l_{r-1} i_{r-1}}
          \left[ \vphantom{\sum_{s=1}^{r-1}}
                 b_{l_{1}j_{\sigma(1)}} \cdots b_{l_{r-1} j_{\sigma(r-1)}}
                 a_{l_{r} i_{r}}  \right.    \nonumber \\[-1mm]
& & \left. \;- \sum_{s=1}^{r-1}   \delta_{l_r l_s}
                b_{l_{1} j_{\sigma(1)}}  \cdots b_{l_{s-1} j_{\sigma(s-1)}}
                b_{l_{s+1} j_{\sigma(s+1)}} \cdots b_{l_{r-1} j_{\sigma(r-1)}}
                h_{i_r j_{\sigma(s)}}
    \right]   b_{l_{r} j_{\sigma(r)}} \\[2mm]
& = & \sum_{\sigma\in{\cal S}_r} \sgn(\sigma)
      \sum_{l_{1},\ldots, l_{r-1} \in [m]}
          a_{l_{1}i_{1}} \cdots a_{l_{r-1} i_{r-1}}
          \left[ \vphantom{\sum_{s=2}^r}
                 b_{l_{1}j_{\sigma(1)}} \cdots b_{l_{r-1} j_{\sigma(r-1)}}
                 (\AT B)_{i_{r}j_{\sigma(r)}}  \right.   \nonumber \\[-1mm]
& & \left. \;- \sum_{s=1}^{r-1} 
             b_{l_{1} j_{\sigma(1)}}  \cdots b_{l_{s-1} j_{\sigma(s-1)}}
             b_{l_{s} j_{\sigma(r)}}
             b_{l_{s+1} j_{\sigma(s+1)}} \cdots b_{l_{r-1} j_{\sigma(r-1)}}
                h_{i_r j_{\sigma(s)}}
    \right]  \\[2mm]
& = & \sum_{\sigma\in{\cal S}_r} \sgn(\sigma)
      \!\!\!
      \sum_{l_{1},\ldots, l_{r-1} \in [m]}
      \!\!
          a_{l_{1}i_{1}} \cdots a_{l_{r-1} i_{r-1}}
          b_{l_{1}j_{\sigma(1)}} \cdots b_{l_{r-1} j_{\sigma(r-1)}}
          \left[ (\AT B)_{i_{r}j_{\sigma(r)}}  +
                 \,\sum_{s=1}^{r-1} h_{i_r j_{\sigma(r)}} \right]
    \nonumber \\[-1mm] \\[2mm]
& = & \sum_{\sigma\in{\cal S}_r} \sgn(\sigma)
      \!\!\!
      \sum_{l_{1},\ldots, l_{r-1} \in [m]}
      \!\!
          a_{l_{1}i_{1}} \cdots a_{l_{r-1} i_{r-1}}
          b_{l_{1}j_{\sigma(1)}} \cdots b_{l_{r-1} j_{\sigma(r-1)}}
      [\AT B + (r-1) h]_{i_{r}j_{\sigma(r)}} 
  \;,
     \nonumber \\[-5mm]
\end{subeqnarray}
where we exchanged $\sigma(s)$ with $\sigma(r)$.
This procedure can be iterated as before to obtain
\begin{subeqnarray}
 \lefteqn{ \sum_{L} (\rowdet \: (\AT)_{IL}) (\det B_{LJ})   }  \\[1mm]
 & = &
 \sum_{\sigma\in{\cal S}_r} \sgn(\sigma)
    \left[\AT B \right]_{i_1 j_{\sigma(1)}}
    \cdots
    \left[\AT B +\,(r-2)\,h \right]_{i_{r-1} j_{\sigma(r-1)}}
    \left[\AT B +\,(r-1)\,h \right]_{i_{r} j_{\sigma(r)}}
        \nonumber\\[-3mm] \\[1mm]
 & = &
 \rowdet  \left[(\AT B)_{I J} + Q_{\rm row} \right]
   \,,
\end{subeqnarray}
which is the desired result of part (b).
\qed

Let us remark that if we care only about the
Capelli identity (i.e., Corollary~\ref{cor.capelli_rectangular} with $h=1$),
then the proof becomes even simpler:
all the discussion about column-pseudo-commutativity is unnecessary
because we have the stronger hypothesis $[a_{ij},a_{kl}]=0$,
so the first steps in the proof proceed exactly as in the commutative case;
and Lemma~\ref{lemma.h_commute} is unnecessary
because the quantities $h_{jl}$ are by hypothesis
central elements of the ring.
The key nontrivial steps in the proof of the Capelli identity
are thus the manipulations of commutators
leading from \reff{eq.mainproof.1} to \reff{eq.mainproof.5}.

\medskip

\begin{example}
   \label{example.itoh}
\rm
Itoh \cite[Theorem~A]{Itoh_04} proves a Capelli-type identity
for some matrices arising from the action of $O(m) \times Sp(n)$
on $m \times n$ matrices ($n$ even).
His result can be obtained as a special case of
Proposition~\ref{prop.capelli_rectangular}${}'$, as follows:
We write $\nu = n/2$ and work in the Weyl algebra $A_{m \times \nu}(K)$
generated by an $m \times \nu$ matrix of indeterminates $X = (x_{ij})$
and the corresponding matrix $\partial = (\partial/\partial x_{ij})$
of differential operators.
Now form the $m \times n$ matrices
$A = (X,\partial)$ and $B = (-\partial, X)$,
which satisfy the commutation relations
$[a_{ij}, b_{kl}] = \delta_{ik} \delta_{jl}$.
Note that $A$ and $B$ are row-commutative
but {\em not}\/ column-pseudo-commutative.
We can therefore apply Proposition~\ref{prop.capelli_rectangular}${}'$(c)
to obtain the identity
\begin{subeqnarray}
   \sum_{\begin{scarray}
            L \subseteq [n] \\
            |L| = r
         \end{scarray}}
   (\det A_{IL}) (\det \: (\BT)_{L J})
   & = &
   \coldet[(A \BT)_{I J} + Q_{\rm col}]   \\[-7mm]
   & = &
   \rowdet[(A \BT)_{I J} + Q_{\rm row}]
   \;,
 \label{eq.prop1.1primec.BIS}
\end{subeqnarray}
in which the matrix $A \BT$ describes the left action of $\ooo(m)$.
But we {\em cannot}\/ obtain a corresponding identity
for the matrix $\AT B$, which describes the right action of $\spsp(n)$,
because $A$ and $B$ fail to be column-pseudo-commutative.
Itoh \cite{Itoh_04} also notes this failure
and interprets it from the point of view of his own proof
(which uses exterior algebra);
he goes on to provide a weaker alternative formula.
\qed
\end{example}

\begin{example}
   \label{example.prop1.1}
{\bf (generalizing Example~\ref{example.easy_noncomm})}
\rm
It is instructive to consider the most general case
with $m=1$, $n=2$ and $r=2$:
here the left-hand side of the identity
is automatically zero (since $r > m$),
but the right-hand side need not vanish
unless we make suitable hypotheses on the matrices
$A = (\alpha,\beta)$ and $B = (\gamma,\delta)$.
We have
$\AT B = \left( \!\! \displaystyle{ \begin{array}{cc}
                                       \alpha\gamma & \alpha\delta \\
                                       \beta\gamma  & \beta\delta
                                     \end{array}
                                  }
           \!\! \right)$,
$H = \left( \!\! \displaystyle{ \begin{array}{cc}
                \gamma\alpha-\alpha\gamma & \delta\alpha-\alpha\delta \\
                \gamma\beta-\beta\gamma   & \delta\beta-\beta\delta
                                     \end{array}
                                  }
           \!\! \right)$
and
$Q_{\rm col} = \left( \!\! \displaystyle{ \begin{array}{cc}
                \gamma\alpha-\alpha\gamma & 0 \\
                \gamma\beta-\beta\gamma   & 0
                                     \end{array}
                                  }
           \!\! \right)$,
so that
$\AT B + Q_{\rm col} = \left( \!\! \displaystyle{ \begin{array}{cc}
                                       \gamma\alpha & \alpha\delta \\
                                       \gamma\beta  & \beta\delta
                                     \end{array}
                                  }
           \!\! \right)$
(a beautiful cancellation!)
and hence
\be
   \coldet(\AT B + Q_{\rm col})
   \;=\;
   \gamma\alpha\beta\delta - \gamma\beta\alpha\delta
   \;=\;
   \gamma \, [\alpha,\beta] \, \delta
   \;.
\ee
In Example~\ref{example.easy_noncomm} we took $\gamma=\delta=1$
and found that the identity fails unless $[\alpha,\beta] = 0$,
i.e.\ $A$ is column-(pseudo-)commutative.
In Itoh's \cite{Itoh_04} $\spsp(n)$ action (specialized to $m=1$ and $n=2$)
we work in the Weyl algebra $A_1(K)$
and have $\alpha=x$, $\beta=d$, $\gamma=-d$ and $\delta=x$,
so that $[\alpha,\beta] = -1 \neq 0$ and the identity again fails.
\qed
\end{example}

\bigskip

{\bf Remarks.}
1.  Part (b) of Proposition~\ref{prop.capelli_rectangular}
is essentially equivalent to part (a)
under a duality that reverses the order of products inside each monomial,
when $R$ is the algebra of noncommutative polynomials (over $\Z$)
in indeterminates $(a_{ij})$ and $(b_{ij})$ with the appropriate relations.
For brevity we omit the details.

2.  Suppose that $[a_{ij}, b_{kl}] = - g_{ik} h_{jl}$
where $G = (g_{ik})$ and $H = (h_{jl})$ are two matrices.
Now make the replacements $A \to PAQ$ and $B \to RAS$,
where $P,Q,R,S$ are matrices whose elements commute
with each other and with those of $A$, $B$, $G$, $H$;
then $G \to PGR^\T$ and $H \to Q^\T HS$.
It follows that Proposition~\ref{prop.capelli_rectangular}
with general $h$ ---
or even the extension to general $g$ and $h$,
{\em provided that $[g_{ik},h_{jl}] = 0$}\/ ---
is not terribly much more general than the traditional case $g=h=\delta$.
However, the case in which $[g_{ik},h_{jl}] \neq 0$
is much more difficult \cite{CS_capelli2}.
\qed

\bigskip

Let us now state and prove a variant of
Proposition~\ref{prop.capelli_rectangular}
in which the hypotheses on the $[a,a]$ or $[b,b]$ commutators are weakened:
instead of assuming that $A$ (or $B$) is column-pseudo-commutative,
we shall assume only that it has column-symmetric commutators.
The price that we have to pay
(as in Lemma~\ref{lemma.coldet.columns}(b))
is that the final identity must be multiplied by a factor of 2 ---
or in other words, the identity holds only up to addition
of an unknown element $x \in R$ satisfying $2x=0$.
Of course, if the ring $R$ is such that $2x=0$ implies $x=0$
(as it is in most applications), then both the hypotheses and the conclusion
are equivalent to those of Proposition~\ref{prop.capelli_rectangular};
but in general rings, the hypotheses are slightly weaker,
as is the conclusion.

\begin{proposition}[noncommutative Cauchy--Binet, variant]
  \label{prop.capelli_rectangular_ANDREA}
Let $R$ be a (not-nec\-es\-sar\-i\-ly-commutative) ring,
and let $A$ and $B$ be $m \times n$ matrices with elements in $R$. 
Suppose that
\be
      [a_{ij}, b_{kl}]   \;=\;  - \delta_{ik} h_{jl}
\ee
where $(h_{jl})_{j,l=1}^n$ are elements of $R$.
Then, for any $I, J \subseteq [n]$ of cardinality $|I|=|J|=r$:
\begin{itemize}
\item[(a)] If $A$ has column-symmetric commutators, then
\be
   2\!\! \sum_{\begin{scarray}
            L \subseteq [m] \\
            |L| = r
         \end{scarray}}
   (\coldet \: (\AT)_{IL}) (\coldet B_{L J})
     \;=\;   2 \coldet[(\AT B)_{I J} + Q_{\rm col}] 
 \label{eq.prop.capelli_rectangular_ANDREA.a}
\ee
where
\be
     (Q_{\rm col})_{\alpha\beta}  \;=\; (r-\beta)\, h_{i_\alpha j_\beta}
\ee
for $1\leq \alpha, \beta \leq r$.
\item[(b)] If $B$ has column-symmetric commutators, then
\be
   2\!\! \sum_{\begin{scarray}
            L \subseteq [m] \\
            |L| = r
         \end{scarray}}
   (\rowdet \: (\AT)_{IL}) (\rowdet B_{L J})
      \;=\;   2 \rowdet[(\AT B)_{I J} + Q_{\rm row}] 
 \label{eq.prop.capelli_rectangular_ANDREA.b}
\ee
where
\be
     (Q_{\rm row})_{\alpha\beta}  \;=\;  (\alpha-1)\, h_{i_\alpha j_\beta}  
\ee
for $1\leq \alpha, \beta \leq r$.
\end{itemize}
%
%
\end{proposition}

In proving Proposition~\ref{prop.capelli_rectangular_ANDREA},
we shall need a variant of Lemma~\ref{lemma.h_commute}:

\begin{lemma}
    \label{lemma.h_commuteweak}
Let $R$ be a (not-necessarily-commutative) ring,
and let $A$ and $B$ be $m \times n$ matrices with elements in $R$.
Suppose that for all $i,k \in [m]$ and $j,l \in [n]$ we have
\begin{subeqnarray}
  [a_{ij}, a_{kl}]   & = &    [a_{il}, a_{kj}]  \\[1mm]
  {}[a_{ij}, b_{kl}] & = &   - \delta_{ik} h_{jl}
 \label{eq.lemma.h_commuteweak.1}
\end{subeqnarray}
where $(h_{jl})_{j,l=1}^n$ are elements of $R$.
Then, for all $i \in [m]$ and $j,l,s \in [n]$ we have
\be
  2\,[a_{ij}, h_{ls}]
    \;=\;
  2\, [a_{il}, h_{js}]
       \;;
 \label{eq.lemma.h_commuteweak.2}
\ee
and if $m \ge 2$ we have
\be
  [a_{ij}, h_{ls}]
    \;=\;
  [a_{il}, h_{js}]
  \;.
 \label{eq.lemma.h_commuteweak.3}
\ee
\end{lemma}

\proof
We consider the Jacobi identity~\reff{Jacobi}
and the corresponding identity with $j$ and $l$ interchanged:
\begin{subeqnarray}
    [a_{ij}, [a_{kl}, b_{rs}]] \,+\,
    [a_{kl}, [b_{rs}, a_{ij}]] \,+\,
    [b_{rs}, [a_{ij}, a_{kl}]]   &=&   0
         \\[1mm]
    {}[a_{il}, [a_{kj}, b_{rs}]] \,+\,
    [a_{kj}, [b_{rs}, a_{il}]] \,+\,
    [b_{rs}, [a_{il}, a_{kj}]]   &=&   0
\end{subeqnarray}
Now subtract the two equations
and use the hypotheses \reff{eq.lemma.h_commuteweak.1}:  we obtain
\be
 \delta_{kr} \left( [a_{ij}, h_{ls}]  -  [a_{il}, h_{js}] \right)  \;=\;
 \delta_{ir}  \left( [a_{kl}, h_{js}] -  [a_{kj}, h_{ls}] \right) \; .
\ee
Now fix arbitrarily the indices $i \in [m]$ and $j,l,s \in [n]$.
Choosing $k=r=i$, we get
\be
  [a_{ij}, h_{ls}] -  [a_{il}, h_{js}]  \;=\;
  [a_{il}, h_{js}] -  [a_{ij}, h_{ls}]
  \;,
\ee
which is \reff{eq.lemma.h_commuteweak.2}.
For $m \ge 2$ we can choose $k=r \in [m]$ different from $i$,
which yields \reff{eq.lemma.h_commuteweak.3}.
\qed

\proofof{Proposition~\ref{prop.capelli_rectangular_ANDREA}}
The proof is almost identical to that of
Proposition~\ref{prop.capelli_rectangular},
with the following changes:
in the argument that $f(l_1,\ldots, l_r) = 0$
whenever two or more arguments take the same value,
we use Lemma~\ref{lemma.coldet.columns}(b)
instead of Lemma~\ref{lemma.coldet.columns}(c);
and in the commutation argument involving $h$,
we use Lemma~\ref{lemma.h_commuteweak} [eq.~\reff{eq.lemma.h_commuteweak.2}]
instead of Lemma~\ref{lemma.h_commute}.
These two weakenings are permissible
since we are multiplying the identity by 2.
\qed

\section{Proof of the Turnbull-type identities}
  \label{sec.symmetric}

In this section we shall prove Propositions~\ref{prop.capelli_symmetric}
and \ref{prop.Turnbull-antisymmetric},
following the same strategy as in the previous proofs.

We begin with a variant of Lemma~\ref{lemma.h_commute}
that will be useful in both the symmetric and antisymmetric cases.

\begin{lemma}
    \label{lemma.h_commute.symmetric}
Let $n \ge 2$, let $R$ be a (not-necessarily-commutative) ring,
and let $A$ and $B$ be $n \times n$ matrices with elements in $R$.
Suppose that for all $i,j,k,l \in [n]$ we have
\begin{subeqnarray}
   [a_{ij}, a_{il}]   & = &    0      \\
   {}[a_{ij}, b_{kl}] & = &   - h_1 \delta_{ik} h_{jl}
                              - h_2 \delta_{il} \delta_{jk}
 \label{eq.lemma.h_commute.symmetric.1}
\end{subeqnarray}
where $h_1$ and $h_2$ are elements of $R$.
Then:
\begin{itemize}
   \item[(a)] If either $i=j$ or $n \ge 3$ (or both),
       we have $[a_{ij}, h_1] = 0$.
   \item[(b)] If $i \neq j$, we have $[a_{ij}, h_1 + h_2] = 0$.
\end{itemize}
\end{lemma}

\proof
For any indices $i,j,k,l,r,s \in [n]$
we have the Jacobi identity
\be
    [a_{ij}, [a_{kl}, b_{rs}]] \,+\,
    [a_{kl}, [b_{rs}, a_{ij}]] \,+\,
    [b_{rs}, [a_{ij}, a_{kl}]]   \;=\;   0
    \;.
 \label{Jacobi.symm}
\ee
Now fix indices $i,j,l$ such that $j \neq l$ (here we use $n \ge 2$)
and take $k=r=i$ and $s=l$.
Using the hypotheses \reff{eq.lemma.h_commute.symmetric.1}, we obtain
\be
   [a_{ij}, h_1 + h_2 \delta_{il}] \;=\; 0
      \quad\hbox{whenever $j \neq l$} \;.
\ee
If $i=j$ or $n \ge 3$, then we can choose 
$l$ to be different from both $i$ and $j$, and conclude that
$[a_{ij}, h_1] = 0$.
If $i \neq j$, then we can choose $l=i$
and conclude that $[a_{ij}, h_1 + h_2] = 0$.
\qed

\proofof{Proposition~\ref{prop.capelli_symmetric}}
Let us first consider the case in which $A$ is symmetric
and we seek a result with $\coldet$.
The result is trivial when $n=1$,
so we can assume that $n \ge 2$
and apply Lemma~\ref{lemma.h_commute.symmetric} with $h_1 = h_2 = h$.
The conclusions of Lemma~\ref{lemma.h_commute.symmetric},
together with the supplementary hypotheses (i) or (ii)
of Proposition~\ref{prop.capelli_symmetric}(a) in case $n=2$,
ensure that $[a_{ij},h] = 0$ for all $i,j$
and hence that we can push all $h$ factors freely to the left.

The first two steps in the proof are identical to those in
Proposition~\ref{prop.easy_noncomm}:  we therefore have
\be
 \sum_{L} (\coldet A_{LI}) \, (\coldet B_{LJ}) 
 \;=\;
 \sum_{\sigma\in{\cal S}_r} \sgn(\sigma) \sum_{l_{1},\cdots, l_{r} \in [m]}
      a_{l_{1} i_{\sigma(1)}} \cdots a_{l_{r} i_{\sigma(r)}}
      b_{l_{1} j_{1}} \cdots b_{l_{r} j_{r}} 
 \;.
\ee
We now push the factor $b_{l_{1} j_{1}}$ to the left:
 \begin{subeqnarray}
\lefteqn{  \sum_{L} (\coldet A_{LI}) \, (\coldet B_{LJ})  }   \\[1mm]
& = &
\sum_{\sigma\in{\cal S}_r} \sgn(\sigma) \sum_{l_{1},\ldots, l_{r} \in [n]}
    a_{l_{1} i_{\sigma(1)}}
    \left[ \vphantom{ \sum_{k=2}^{m}}
           b_{l_{1} j_1} a_{l_{2} i_{\sigma(2)}} \cdots a_{l_{r} i_{\sigma(r)}}
    \right.  \nonumber \\
& & \;- h \, \left. \sum_{s=2}^{r} (\delta_{i_{\sigma(s)} j_1} \delta_{l_1 l_s}
                              + \delta_{i_{\sigma(s)} l_{1}} \delta_{l_s j_1} )
           \, a_{l_{2} i_{\sigma(2)}} \cdots a_{l_{s-1} i_{\sigma(s-1)}}
              a_{l_{s+1} i_{\sigma(s+1)}} \cdots a_{l_{r} i_{\sigma(r)}}
    \right] \times   \nonumber \\
& & 
     \qquad\qquad  b_{l_2 j_2} \cdots b_{l_r j_r}   \\[2mm]
& = &
 \sum_{\sigma\in{\cal S}_r} \sgn(\sigma)  \sum_{l_{2},\ldots, l_{r} \in [n]}
    \left[ \vphantom{ \sum_{k=2}^{m}}
           (\AT B)_{i_{\sigma(1)} j_1}
           a_{l_{2} i_{\sigma(2)}} \cdots a_{l_{r} i_{\sigma(r)}} 
    \right.  \nonumber \\
&&  \;- h \sum_{s=2}^{r} \delta_{i_{\sigma(s)} j_1}
           a_{l_{2} i_{\sigma(2)}} \cdots a_{l_{s-1} i_{\sigma(s-1)}}
           a_{l_{s} i_{\sigma(1)}}
           a_{l_{s+1} i_{\sigma(s+1)}} \cdots a_{l_{r} i_{\sigma(r)}}
      \nonumber \\
&&  \;- h \left. \sum_{s=2}^{r} \delta_{l_s j_1}
           a_{l_{2} i_{\sigma(2)}} \cdots a_{l_{s-1} i_{\sigma(s-1)}}
           a_{i_{\sigma(s)} i_{\sigma(1)}}
           a_{l_{s+1} i_{\sigma(s+1)}} \cdots a_{l_{r} i_{\sigma(r)}}
     \right]  b_{l_2 j_2} \cdots b_{l_r j_r}   \quad \\[2mm]
&=&
 \sum_{\sigma\in{\cal S}_r} \sgn(\sigma)
       \left[ (\AT  B)_{i_{\sigma(1)} j_1}  +
               h\,\sum_{s=2}^{r} \delta_{i_{\sigma(1)} j_1} \right] 
       \sum_{l_{2},\ldots, l_{r} \in [n]}
            a_{l_{2} i_{\sigma(2)}} \cdots a_{l_{r} i_{\sigma(r)}}
            b_{l_{2} j_{2}} \cdots b_{l_{r} j_{r}}
     \nonumber \\[-2mm]   \\[2mm]
&=&
 \sum_{\sigma\in{\cal S}_r} \sgn(\sigma)
       \left[\AT B  + h(r-1) \right]_{i_{\sigma(1)} j_1}
       \sum_{l_{2},\ldots, l_{r} \in [n]}
            a_{l_{2} i_{\sigma(2)}} \cdots a_{l_{r} i_{\sigma(r)}}
            b_{l_{2} j_{2}} \cdots b_{l_{r} j_{r}}
     \;,
\end{subeqnarray}
where in the interchange of $\sigma(1)$ with $\sigma(s)$
we dropped out the third summand within square brackets
because it is symmetric, hence vanishing when summed over $\sigma$
with a factor $\sgn(\sigma)$.
Iterating this procedure, we obtain the desired final result
\begin{subeqnarray}
 \lefteqn{ \sum_{L} (\det  A_{LI}) (\det B_{LJ})   }  \\[1mm]
 & = &
 \sum_{\sigma\in{\cal S}_r} \sgn(\sigma)
    \left[\AT B  + \,(r-1)\,h \right]_{i_{\sigma(1)}j_{1}}
    \left[\AT B  + \,(r-2)\,h \right]_{i_{\sigma(2)}j_{2}}
    \cdots
    \left[\AT B \right]_{i_{\sigma(r)}j_{r}}     \nonumber\\[-4mm] \\[1mm]
 & = &
 \coldet  \left[(\AT B)_{I J} + Q_{\rm col} \right]
   \,.
\end{subeqnarray}

The proof of part (b) is similar.
\qed

{\bf Remark.}
In the proof of Proposition~\ref{prop.capelli_symmetric}
contained in the first preprint version of this paper,
we inadvertently assumed without justification that $[a_{ij},h] = 0$.\footnote{
   We realized our error when we read \cite{Chervov_09}
   and in particular their Lemmas~21 and 25,
   which correspond to our Lemma~\ref{lemma.h_commute.symmetric}
   albeit restricted to the case where $2x=0$ implies $x=0$ in the ring $R$.
}
We have now repaired the error by including
Lemma~\ref{lemma.h_commute.symmetric}
and adding the extra conditions (i) or (ii)
in Proposition~\ref{prop.capelli_symmetric} when $n=2$.
The importance of these extra conditions
is illustrated by the following example.
\qed

\begin{example}
   \label{example.2by2.symmetric}
\rm
Let us consider the general case of $2 \times 2$ matrices (i.e.\ $n=2$)
under the hypotheses that $A$ is column-pseudo-commutative
and that $[a_{ij}, b_{kl}] = -  h (\delta_{ik} \delta_{jl} +
                                   \delta_{il} \delta_{jk} )$.
With $Q = \left( \!\! \displaystyle{ \begin{array}{cc}
                                          h & 0 \\
                                          0 & 0
                                      \end{array}
                                    }
           \!\! \right)$,
we have
\be
  \coldet(\AT B + Q) \,-\, (\coldet \AT )(\coldet B)
  \;=\;
     - [a_{12}, h] b_{12} \,+\, (a_{21} - a_{12}) h b_{12}
  \;.
 \label{eq.ab22.symmetric}
\ee
The second term illustrates the importance of assuming that $A$ is symmetric. 
The first term illustrates the importance of assuming that
$[a_{12}, h] = 0$ when $n=2$
(i.e.\ when this is not guaranteed by Lemma~\ref{lemma.h_commute.symmetric}).

To see that $[a_{12}, h] \neq 0$
and $\hbox{\reff{eq.ab22.symmetric}} \neq 0$
can actually arise,
let $R$ be the ring of $2 \times 2$ matrices with elements in
the field $GF(2)$, and consider the elements
 $\alpha = \left( \!\! \displaystyle{ \begin{array}{cc}
                                          1 & 0 \\
                                          0 & 0
                                      \end{array}
                                    }
           \!\! \right)$
 and
 $\beta  = \left( \!\! \displaystyle{ \begin{array}{cc}
                                          0 & 1 \\
                                          1 & 0
                                      \end{array}
                                    }
           \!\! \right)$
in $R$.
Note that $[\alpha,\beta] = \beta$.
Now let
$A =  \left( \!\! \displaystyle{ \begin{array}{cc}
                                     0 & \alpha \\
                                     \alpha  & 0
                                 \end{array}
                               }
      \!\! \right)$
and
$B =  \left( \!\! \displaystyle{ \begin{array}{cc}
                                     0 & \beta \\
                                     \beta  & 0
                                 \end{array}
                               }
      \!\! \right)$.
These matrices satisfy $[a_{ij},a_{kl}] = [b_{ij},b_{kl}] = 0$
and $[a_{ij}, b_{kl}] = -  h (\delta_{ik} \delta_{jl} +
                              \delta_{il} \delta_{jk} )$
with $h=\beta$ (this works when $ij=kl=11$ or $22$ because $2h=0$).
But $[a_{12},h] = [\alpha,\beta] = \beta \neq 0$
and $\hbox{\reff{eq.ab22.symmetric}} = \beta^2 \neq 0$.
\qed
\end{example}

\proofof{Proposition~\ref{prop.Turnbull-antisymmetric}}
Let us first consider the case in which $A$ is antisymmetric off-diagonal
and we seek a result with $\colper$.
We follow closely the proof of the previous
Proposition~\ref{prop.capelli_symmetric},
starting from the expression in which the $a$'s and $b$'s are ordered
and first pushing the factor $b_{l_{1} j_{1}}$ to the left:
\begin{subeqnarray}
\lefteqn{  \sum_{\sigma\in{\cal S}_r}  \sum_{l_{1},\cdots, l_{r} \in [n]}
      a_{l_{1} i_{\sigma(1)}} \cdots a_{l_{r} i_{\sigma(r)}}
      b_{l_{1} j_{1}} \cdots b_{l_{r} j_{r}}  }   \\[1mm]
& = &
\sum_{\sigma\in{\cal S}_r}  \sum_{l_{1},\ldots, l_{r} \in [n]}
    a_{l_{1} i_{\sigma(1)}}
    \left[ \vphantom{ \sum_{k=2}^{m}}
           b_{l_{1} j_1} a_{l_{2} i_{\sigma(2)}} \cdots a_{l_{r} i_{\sigma(r)}}
    \right.  \nonumber \\
& & \;- h \, \left. \sum_{s=2}^{r} (\delta_{i_{\sigma(s)} j_1} \delta_{l_1 l_s}
                              - \delta_{i_{\sigma(s)} l_{1}} \delta_{l_s j_1} )
           \, a_{l_{2} i_{\sigma(2)}} \cdots a_{l_{s-1} i_{\sigma(s-1)}}
              a_{l_{s+1} i_{\sigma(s+1)}} \cdots a_{l_{r} i_{\sigma(r)}}
    \right] \times   \nonumber \\
& & 
     \qquad\qquad  b_{l_2 j_2} \cdots b_{l_r j_r}   \\[2mm]
& = &
 \sum_{\sigma\in{\cal S}_r}  \sum_{l_{2},\ldots, l_{r} \in [n]}
    \left[ \vphantom{ \sum_{k=2}^{m}}
           (\AT B)_{i_{\sigma(1)} j_1}
           a_{l_{2} i_{\sigma(2)}} \cdots a_{l_{r} i_{\sigma(r)}} 
    \right.  \nonumber \\
&&  \;- h \sum_{s=2}^{r} \delta_{i_{\sigma(s)} j_1}
           a_{l_{2} i_{\sigma(2)}} \cdots a_{l_{s-1} i_{\sigma(s-1)}}
           a_{l_{s} i_{\sigma(1)}}
           a_{l_{s+1} i_{\sigma(s+1)}} \cdots a_{l_{r} i_{\sigma(r)}}
      \nonumber \\
&&  \;+ h \left. \sum_{s=2}^{r} \delta_{l_s j_1}
           a_{l_{2} i_{\sigma(2)}} \cdots a_{l_{s-1} i_{\sigma(s-1)}}
           a_{i_{\sigma(s)} i_{\sigma(1)}}
           a_{l_{s+1} i_{\sigma(s+1)}} \cdots a_{l_{r} i_{\sigma(r)}}
     \right]  b_{l_2 j_2} \cdots b_{l_r j_r}   \quad \\[2mm]
&=&
 \sum_{\sigma\in{\cal S}_r} 
       \left[ (\AT  B)_{i_{\sigma(1)} j_1}  -
               h\,\sum_{s=2}^{r} \delta_{i_{\sigma(1)} j_1} \right] 
       \sum_{l_{2},\ldots, l_{r} \in [n]}
            a_{l_{2} i_{\sigma(2)}} \cdots a_{l_{r} i_{\sigma(r)}}
            b_{l_{2} j_{2}} \cdots b_{l_{r} j_{r}}
     \nonumber \\[-2mm]   \\[2mm]
&=&
 \sum_{\sigma\in{\cal S}_r} 
       \left[\AT B  - h(r-1) \right]_{i_{\sigma(1)} j_1}
       \sum_{l_{2},\ldots, l_{r} \in [n]}
            a_{l_{2} i_{\sigma(2)}} \cdots a_{l_{r} i_{\sigma(r)}}
            b_{l_{2} j_{2}} \cdots b_{l_{r} j_{r}}
     \;,
\end{subeqnarray}
where now in the interchange of $\sigma(1)$ with $\sigma(s)$
there is {\em no}\/ change of sign in the second summand
because there is {\em no}\/ $\sgn(\sigma)$ factor,
and we dropped out the third summand within square brackets
because it is {\em anti}\/symmetric
and hence vanishing when summed over $\sigma$
(note that we used here only the {\em off-diagonal}\/ antisymmetry of $A$).
Iterating this procedure, we obtain the desired final result
\begin{subeqnarray}
\lefteqn{  \sum_{\sigma\in{\cal S}_r}  \sum_{l_{1},\cdots, l_{r} \in [n]}
      a_{l_{1} i_{\sigma(1)}} \cdots a_{l_{r} i_{\sigma(r)}}
      b_{l_{1} j_{1}} \cdots b_{l_{r} j_{r}}  }   \\[1mm]
 & = &
 \sum_{\sigma\in{\cal S}_r} 
    \left[\AT B  - \,(r-1)\,h \right]_{i_{\sigma(1)}j_{1}}
    \left[\AT B  - \,(r-2)\,h \right]_{i_{\sigma(2)}j_{2}}
    \cdots
    \left[\AT B \right]_{i_{\sigma(r)}j_{r}}     \nonumber\\[-4mm] \\[1mm]
 & = &
 \colper \left[(\AT B)_{I J} - Q_{\rm col} \right]
   \,.
\end{subeqnarray}

The proof of part (b) is similar.
\qed

\medskip
\noindent
{\bf Remark.}
In the proof of Proposition~\ref{prop.Turnbull-antisymmetric}
contained in the first preprint version of this paper,
we inadvertently assumed without justification that $[a_{ij},h] = 0$
or $[b_{ij},h] = 0$.  We have now repaired the error
by including these hypotheses explicitly.
%
Alternatively, if we assume that $[a_{ij},a_{il}] = 0$
for all $i,j,l$, then Lemma~\ref{lemma.h_commute.symmetric}
with $h_1 = -h_2 = h$ implies that $[a_{ij},h]=0$ for all $i,j$
whenever $n \ge 3$ (and analogously for $B$).
\qed

\section{Further generalizations?}
  \label{sec.further}

In this section we investigate whether any of our results
can be further generalized.
We always assume, for simplicity, that
$[a_{ij},a_{kl}]= 0$ and $[b_{ij},b_{kl}]= 0$ for all $i,j,k,l$,
and we ask the following question: 
\begin{quote}
\noindent
{\em In which cases does there exist a diagonal matrix $Q$ such that
\be
    (\det A)(\det B)  \;=\; \coldet( A^\T B + Q)
    \,\hbox{\it ?}
 \label{eq.further.1}
\ee
And what are the elements of $Q$?}
\end{quote}
More specifically, we investigate two types of possible extensions:
\begin{itemize}
   \item  Can we allow commutators
    $[a_{ij},b_{kl}] =
       -h_1 \delta_{ik} \delta_{jl} - h_2 \delta_{il} \delta_{jk}$
    beyond the ordinary case ($h_2=0$) and the symmetric case ($h_2=h_1$)?
   \item  In the symmetric case $h_2=h_1$, is it really necessary
    that the matrix $A$ be symmetric?
\end{itemize}
Our approach is to perform exact calculations for small matrices,
in order to conjecture possible identities and to rule out others.
We shall find, not surprisingly, that there exists an antisymmetric case
($h_2=-h_1$) corresponding to the Howe--Umeda--Kostant--Sahi identity
\cite{Howe_91,Kostant_91}
with $n$ even;  but it appears that there exists also {\em another}\/
antisymmetric case, with a {\em different}\/ ``quantum correction'',
when $n$ is odd.

We shall perform our calculations in two alternative frameworks:

\medskip

{\bf Abstract framework.}
We work in the ring $\Z\<A,B\>[h_1,h_2]/\scrr$
generated by noncommuting indeterminates
$A = (a_{ij})_{i,j=1}^n$ and $B = (b_{ij})_{i,j=1}^n$
and commuting indeterminates $h_1$ and $h_2$
modulo the two-sided ideal $\scrr$ generated by the
relations $[a_{ij},a_{kl}]= 0$, $[b_{ij},b_{kl}]= 0$ and
$[a_{ij},b_{kl}] =
   -h_1 \delta_{ik} \delta_{jl} - h_2 \delta_{il} \delta_{jk}$.
We then introduce a matrix $Q = \diag(q_1,\ldots,q_n)$
of central elements and expand out the polynomial
\be
   f(A,B,Q)  \;\bydef\;
   (\det A)(\det B)  \,-\, \coldet( A^\T B + Q )
   \;.
\ee
We ask whether $q_1,\ldots,q_n$ can be chosen so that either
\begin{itemize}
   \item $f(A,B,Q) \equiv 0$, or
   \item $f(A,B,Q) \equiv 0$ modulo the ideal corresponding to the
      symmetry or antisymmetry of $A$ and/or $B$.
      More precisely, there are six nontrivial ideals to be considered:
      \begin{itemize}
           \item $A$ symmetric
           \item $B$ symmetric
           \item $A$ and $B$ symmetric
           \item $A$ antisymmetric
           \item $B$ antisymmetric
           \item $A$ and $B$ antisymmetric
      \end{itemize}
      In the first three cases we have $h_2 = h_1$,
      while in the latter three cases we have $h_2 = -h_1$.
      (If $A$ and $B$ have opposite symmetries,
       then $h_1 = h_2 = 0$ and we are back in the commutative case.)
\end{itemize}

\medskip

{\bf Concrete framework.}
We work in the Weyl algebra $A_{n \times n}(K)$
over a field $K$ of characteristic 0 (e.g.\ $\Q$, $\R$ or $\C$)
generated by a matrix $X=(x_{ij})_{i,j=1}^n$ of indeterminates
and the corresponding matrix $\partial = (\partial/\partial x_{ij})_{i,j=1}^n$
of partial derivatives.
The matrices $X$ and $Y = h\partial$ (where $h \in K$, $h \neq 0$)
satisfy the commutation relations
\be
[x_{ij},x_{kl}]  = 0\, ; \qquad
[y_{ij},y_{kl}]  = 0 \, ; \qquad
[x_{ij},y_{kl}]  = - h\, \delta_{ik}\,\delta_{jl}   \;.
   \label{weyl} 
\ee
Now consider the matrices $A= (a_{ij})_{i,j=1}^n$ and $B= (b_{ij})_{i,j=1}^n$
defined by
\be
A  = \alpha X + \beta X^\T \, ; \qquad
B  = \gamma Y + \delta Y^\T
\ee
(where $\alpha,\beta,\gamma,\delta \in K$), so that
\begin{subeqnarray}
  [a_{ij},a_{kl}]   & = &   0   \\
  {}[b_{ij},b_{kl}] & = &   0   \\
  {}[a_{ij},b_{kl}] & = &
    - h \left[ (\alpha \gamma + \beta \delta) \, \delta_{ik} \delta_{jl}
         \,+\, (\alpha \delta + \beta \gamma) \, \delta_{il} \delta_{jk}
        \right]
  \,.
 \label{gen}
\end{subeqnarray}
Some cases of special interest are:
\begin{itemize}
\item $\alpha \delta + \beta \gamma = 0$, so that \reff{gen} reduces
    to the starting case~\reff{weyl} with the replacement
    $h \to h \,(\alpha \gamma + \beta \delta)$;
\item $\beta = \alpha$, so that $A$ is symmetric
    and $\alpha \gamma + \beta \delta=\alpha \delta + \beta \gamma
         =\alpha(\gamma+\delta)$;
\item $\delta = \gamma$, so that $B$ is symmetric
    and $\alpha \gamma + \beta \delta=\alpha \delta + \beta \gamma
         =(\alpha + \beta)\gamma$;
\item $\beta = \alpha$ and $\delta = \gamma$, so that both $A$ and $B$
    are symmetric and
    $\alpha \gamma + \beta \delta=\alpha \delta + \beta \gamma= 2 \alpha\gamma$;
\item $\beta = -\alpha$, so that $A$ is antisymmetric
    and $\alpha \gamma + \beta \delta=- (\alpha \delta + \beta \gamma)
    =\alpha(\gamma-\delta)$;
\item $\delta = -\gamma$, so that $B$ is antisymmetric
    and $\alpha \gamma + \beta \delta=-(\alpha \delta + \beta \gamma)
    =(\alpha - \beta)\gamma$;
\item $\beta = -\alpha$ and $\delta = -\gamma$, so that  both $A$ and $B$
    are antisymmetric and $\alpha \gamma + \beta \delta=
     -(\alpha \delta + \beta \gamma)= 2 \alpha \gamma$.
\end{itemize}
%

\medskip

We shall use the abstract framework to conjecture possible identities,
and the concrete framework to find counterexamples.

In the abstract framework we know the following general solutions:
\begin{itemize}
   \item[(a)]  Capelli:  For any $n$, with $h_2=0$ and $q_i = (n-i) h_1$.
   \item[(b)]  Turnbull:  For any $n$, with $h_2=h_1$, $A$ symmetric
       and $q_i = (n-i) h_1$.
       [Turnbull requires symmetry of both $A$ and $B$,
        but we have shown in Proposition~\ref{prop.capelli_symmetric}
        that only the symmetry of $A$ is needed.]
   \item[(c)]  Howe--Umeda--Konstant--Sahi:
       For any {\em even}\/ $n$, with $h_2 = -h_1$, $A$ and $B$ antisymmetric
       and $q_i = (n-i-1) h_1$.
       [But we suspect that only $B$ need be assumed antisymmetric: see below.]
\end{itemize}
Our results lead us to {\em conjecture}\/ a fourth general solution:
\begin{itemize}
   \item[(d)]  For any {\em odd}\/ $n$, with $h_2 = -h_1$,
       $A$ or $B$ antisymmetric (or both) and $q_i = (n-i)h_1$.
\end{itemize}
We shall show
--- with the aid of the symbolic-algebra package {\sc Mathematica} ---
that for $n \le 5$ there are no other solutions,
except for a special antisymmetric solution at $n=2$.

\bigskip

{\bf Case \boldmath $n=1$.}  This case is trivial and we have $Q=0$.

\bigskip

{\bf Case \boldmath $n=2$.}
For $n=2$ we find the following solutions:
\begin{itemize}
   \item Capelli: $h_2 = 0$, $q_1 = h_1$, $q_2 = 0$.
   \item Turnbull: $h_2 = h_1$, $A$ symmetric, $q_1 = h_1$, $q_2 = 0$.
   \item Howe--Umeda--Konstant--Sahi: $h_2 = -h_1$,
              $B$ antisymmetric, $q_1 = 0$, $q_2 = -h_1$.
   \item New antisymmetric solution:  $h_2 = -h_1$,
              $A$ antisymmetric, $q_1 = h_1$, $q_2 = 0$.
\end{itemize}
We see, therefore, that for $n=2$ the Howe--Umeda--Konstant--Sahi solution
requires only the antisymmetry of $B$, not of $A$.
Furthermore, there is an additional antisymmetric solution
with a different quantum correction;
this solution does not, however, seem to generalize to larger even $n$
(see $n=4$ below).

There are no solutions for $n=2$ besided the ones listed above.
Indeed, even in the more restrictive concrete framework
we can show that unless we are in one of the foregoing cases,
there is no choice of $q_1$ and $q_2$ for which \reff{eq.further.1} holds.

\bigskip

{\bf Case \boldmath $n=3$.}
For $n=3$ we find the following solutions:
\begin{itemize}
   \item Capelli: $h_2 = 0$, $q_i = (n-i) h_1$.
   \item Turnbull: $h_2 = h_1$, $A$ symmetric, $q_i = (n-i) h_1$.
   \item New antisymmetric solution:  $h_2 = -h_1$,
              either $A$ or $B$ (or both) antisymmetric, $q_i = (n-i) h_1$.
\end{itemize}

There are no solutions for $n=3$ besided the ones listed above.
Indeed, even in the more restrictive concrete framework
we can show that unless we are in one of the foregoing cases,
there is no choice of $q_1$, $q_2$ and $q_3$
for which \reff{eq.further.1} holds.

\bigskip

{\bf Case \boldmath $n=4$.}
For $n=4$ we find the following solutions:
\begin{itemize}
   \item Capelli: $h_2 = 0$, $q_i = (n-i) h_1$.
   \item Turnbull: $h_2 = h_1$, $A$ symmetric, $q_i = (n-i) h_1$.
   \item Howe--Umeda--Konstant--Sahi: $h_2 = -h_1$,
              $B$ antisymmetric, $q_i = (n-i-1) h_1$.
\end{itemize}

There are no solutions for $n=4$ besided the ones listed above.
(With our available computer facilities
we were able to perform the computation for $n=4,5$
only in the abstract framework, not in the concrete framework.)

\bigskip

{\bf Case \boldmath $n=5$.}
For $n=5$ we find the following solutions:
\begin{itemize}
   \item Capelli: $h_2 = 0$, $q_i = (n-i) h_1$.
   \item Turnbull: $h_2 = h_1$, $A$ symmetric, $q_i = (n-i) h_1$.
   \item New antisymmetric solution:  $h_2 = -h_1$,
              either $A$ or $B$ (or both) antisymmetric, $q_i = (n-i) h_1$.
\end{itemize}

There are no solutions for $n=5$ besided the ones listed above.

\bigskip

These results for $n \le 5$ lead us to make the following conjectures:

\begin{conjecture}[generalized Howe--Umeda--Konstant--Sahi]
   \label{conj1}
If $n$ is {\em even}\/ and $B$ is antisymmetric (but $A$ can be arbitrary),
then \reff{eq.further.1} holds with
\be
   q_i \;=\; (n-i-1) h_1  \;.
\ee
\end{conjecture}

\begin{conjecture}[new antisymmetric solution]
   \label{conj2}
If $n$ is {\em odd}\/ and at least one of the matrices $A$ and $B$ is
antisymmetric, then \reff{eq.further.1} holds with
\be
   q_i \;=\; (n-i) h_1  \;.
\ee
\end{conjecture}

\clearpage 

\appendix
\section{A generalized ``Cayley'' identity}

In this appendix we use Proposition~\ref{prop.capelli_rectangular}
to prove a generalization of the ``Cayley'' identity \reff{eq.intro.1}.

\begin{proposition}[generalized Cayley identity]
   \label{prop.genCayley}
Let $R$ be a (not-necessarily-com\-mu\-ta\-tive) ring,
and let $A$ and $B$ be $n \times n$ matrices with elements in $R$.
Suppose that
\begin{subeqnarray}
       [a_{ij}, a_{kl}]   & = &  0  \\[1mm]
       [a_{ij}, b_{kl}]   & = &  - \delta_{ik} h_{jl}  \\[1mm]
       [a_{ij}, h_{kl}]   & = &  0
\end{subeqnarray}
where $H = (h_{jl})_{j,l=1}^n$ is a matrix with elements in $R$.
Then, for any $I, J \subseteq [n]$ of cardinality $|I|=|J|=r$
and any nonnegative integer $s$, we have
\begin{eqnarray}
    \sum_{\begin{scarray}
            L \subseteq [n] \\
            |L| = r
         \end{scarray}}
   (\det \: (\AT)_{IL}) (\coldet B_{L J})  (\det A)^s
     &=&   (\det A)^s\, \coldet[(\AT B + sH)_{I J} + Q_{\rm col}]
   \;, \!\!\!\!
   \nonumber \\[-7mm]
  \label{eq.prop.genCayley}
\end{eqnarray}
where
\be
     (Q_{\rm col})_{\alpha\beta}  \;=\; (r-\beta)\, h_{i_\alpha j_\beta}
\ee
for $1\leq \alpha, \beta \leq r$.
\end{proposition}

The identity \reff{eq.prop.genCayley} clearly generalizes
the identity \reff{eq.prop.capelli_rectangular.a}
from Proposition~\ref{prop.capelli_rectangular}(a),
to which it reduces when $s=0$
(and $A$ is commutative: see also Remarks 2 and 3 below).

The key fact needed in the proof of Proposition~\ref{prop.genCayley}
is the following identity:

\begin{lemma}
   \label{lemma.genCayley}
Let $A$ and $B$ be as in Proposition~\ref{prop.genCayley}.
Then, for all $i,j \in [n]$ and all nonnegative integers $s$, we have
\be
  [(\AT B)_{ij}, (\det A)^s]
  \;=\;
  s \, h_{ij} \, (\det A)^s
  \;=\;
  s \, (\det A)^s \, h_{ij} 
  \;.
 \label{eq.lemma.genCayley}
\ee
\end{lemma}   

\proof
A simple computation using the hypotheses $[a,a]=0$ shows that
\be
   [(\AT B)_{ij}, a_{kl}]  \;=\;  a_{ki} h_{lj}
   \;.
 \label{eq.proof.lemma.genCayley.1}
\ee
We therefore have (using the hypotheses $[a,h]=0$)
\begin{subeqnarray}
  & & [(\AT B)_{ij}, \det A]
      \nonumber \\
  & & \quad=\;
  \sum_{\sigma \in \scrs_n} \sum_{r=1}^n \sgn(\sigma) \, a_{\sigma(1) 1}  
\cdots a_{\sigma(r-1), r-1} \ [ (\AT B)_{ij}, a_{\sigma(r) r}] \,
    a_{\sigma(r+1), r+1} \cdots a_{\sigma(n) n} 
         \qquad\quad \\
  & & \quad=\;
  \sum_{r=1}^n h_{rj}\, \sum_{\sigma \in \scrs_n}
       \sgn(\sigma) \, a_{\sigma(1) 1} \cdots
      a_{\sigma(r-1), r-1} \, a_{\sigma(r) i} \, a_{\sigma(r+1), r+1} \cdots
      a_{\sigma(n) n} \\
  & & \quad=\;
   \sum_{r=1}^n h_{rj}\, \delta_{ir} \, (\det A) \\
  & & \quad=\;
   h_{ij} \, (\det A)
   \;,
 \label{eq.proof.lemma.genCayley.2}
\end{subeqnarray}
where the third equality used the fact that the determinant
of a (commutative) matrix with two equal columns is zero
(so that the terms with $i \neq r$ vanish).
This proves \reff{eq.lemma.genCayley} for the case $s=1$.
The general case easily follows by induction using
$[x,yz] = [x,y] z + y[x,z]$ along with the hypotheses $[a,h]=0$.
\qed

\proofof{Proposition~\ref{prop.genCayley}}
By Proposition~\ref{prop.capelli_rectangular} we have
\be
    \sum_{\begin{scarray}
            L \subseteq [n] \\
            |L| = r
         \end{scarray}}
   (\det \: (\AT)_{IL}) (\coldet B_{L J})  (\det A)^s
     \;=\;   \coldet[(\AT B)_{I J} + Q_{\rm col}] \, (\det A)^s
   \;.
\ee
Now let us work on the right-hand side.
By Lemma~\ref{lemma.genCayley} we have
\be
   [((\AT B)_{I J} + Q_{\rm col})_{\alpha\beta}, (\det A)^s]
   \;=\;
   s \, (\det A)^s \, h_{i_\alpha j_\beta} 
\ee
(since the matrix elements of $Q_{\rm col}$ commute with those of $A$),
or in other words
\be
   ((\AT B)_{I J} + Q_{\rm col})_{\alpha\beta} \, (\det A)^s
   \;=\;
   (\det A)^s \, ((\AT B + sH)_{I J} + Q_{\rm col})_{\alpha\beta}
   \;.
 \label{eq.cayley.commutator} 
\ee
Now expand out $\coldet[(\AT B)_{I J} + Q_{\rm col}]$
and right-multiply it by $(\det A)^s$;
repeatedly using \reff{eq.cayley.commutator} to push $(\det A)^s$
to the left, we obtain
[using the shorthands $M = (\AT B)_{I J} + Q_{\rm col}$
 and $M' = (\AT B + sH)_{I J} + Q_{\rm col}$]
\begin{subeqnarray}
   (\coldet M) \, (\det A)^s
   & = &
   \sum_{\sigma \in \scrs_r} \sgn(\sigma) \,
      M_{\sigma(1) 1} \cdots M_{\sigma(r) r} \, (\det A)^s
         \qquad\quad \\[1mm]
   & = &
   \sum_{\sigma \in \scrs_r} \sgn(\sigma) \,
      (\det A)^s \, M'_{\sigma(1) 1} \cdots M'_{\sigma(r) r}
         \\[1mm]
   & = &
   (\det A)^s \, (\coldet M')
   \;,
\end{subeqnarray}
as was to be proved.
\qed

Specializing Proposition~\ref{prop.genCayley}
to the Weyl algebra $R = A_{n \times n}(K)$
over a field $K$ of characteristic 0, we obtain:

\begin{corollary}[Cayley identity, Weyl algebra version]
  \label{cor.genCayley.1}
$\quad\;$
Let $X = (x_{ij})_{i,j=1}^n$ be a square matrix of commuting  
indeterminates, and let $\partial = (\partial/\partial x_{ij})_{i,j=1}^n$  
be the corresponding matrix of partial derivatives.
Then, for any $I, J \subseteq [n]$ of cardinality $|I|=|J|=r$
and any nonnegative integer $s$, we have
\be
    \sum_{\begin{scarray}
            L \subseteq [n] \\
            |L| = r
         \end{scarray}}
   (\det \: (\XT)_{IL}) (\det \partial_{L J})  (\det X)^s
     \;=\;   (\det X)^s\, \coldet[(\XT \partial)_{I J} + Q_{\rm col}(s)]
   \;,
 \label{eq.cor.genCayley.1a}
\ee
where
\be
     Q_{\rm col}(s)_{\alpha\beta}  \;=\;
        (s+r-\beta) \, \delta_{i_\alpha j_\beta}
\ee
for $1\leq \alpha, \beta \leq r$.
In particular, for $I=J=[n]$ we have
\be
    (\det X) (\det \partial) (\det X)^s
      \;=\;   (\det X)^s\, \coldet(\XT \partial + \Delta)
 \label{eq.cor.genCayley.1b}
\ee
where
\be
\Delta_{ij} = (s+n-j)\, \delta_{ij}\, .
\ee
\end{corollary}

Finally, let us take the identity \reff{eq.cor.genCayley.1b}
in the Weyl algebra and apply it to the constant polynomial $1 \in K[X]$;
then all the derivatives $\partial_{ij}$ annihilate 1,
and by removing an overall factor $(\det X)$ from the left
we obtain the usual Cayley identity:

\begin{corollary}[Cayley identity, traditional version]
  \label{cor.genCayley.2}
Let $X = (x_{ij})_{i,j=1}^n$ be a square matrix of commuting
indeterminates, and let $\partial = (\partial/\partial x_{ij})_{i,j=1}^n$
be the corresponding matrix of partial derivatives.
Then, for any nonnegative integer $s$, we have
\be
    (\det \partial) (\det X)^s
      \;=\; s (s+1) \dots (s+n-1)\, (\det X)^{s-1} \,.
 \label{eq.cor.genCayley.2}
\ee
\end{corollary}

{\bf Remarks.}
1.  Let us stress that \reff{eq.cor.genCayley.1a}/\reff{eq.cor.genCayley.1b}
are identities in the Weyl algebra;
they can be applied to {\em any}\/ polynomial $P \in K[X]$
to obtain an identity in the polynomial algebra.
They are therefore stronger than
the traditional Cayley identity \reff{eq.cor.genCayley.2},
which corresponds to taking $P=1$ only.
This example suggests that it might be fruitful
to investigate more generally whether identities of
Bernstein--Sato type \cite{Bernstein_72,Bjork_79,Coutinho_95,Krause_00}
can be extended in a useful way to identities in the Weyl algebra.

2. The hypothesis in Proposition~\ref{prop.genCayley}
that $[a_{ij},h_{kl}]=0$ can be avoided when $n \ge 2$,
since by a simple modification of Lemma~\ref{lemma.h_commute}
it can be shown that $[a_{ij},a_{kl}]=0$
and $[a_{ij}, b_{kl}] = -\delta_{ik} h_{jl}$ for all $i,j,k,l$
{\em implies}\/ $[a_{ij},h_{kl}]=0$ for all $i,j,k,l$,
{\em provided that $n \ge 2$}\/.
For $n=1$, by contrast, this implication is in general false,
and neither Proposition~\ref{prop.genCayley} nor Lemma~\ref{lemma.genCayley}
holds when $[a,h] \neq 0$.

3. It would be interesting to know whether the hypothesis
on $[a,a]$ commutators in Proposition~\ref{prop.genCayley}
can be weakened from commutativity to column-pseudo-commutativity
(of course replacing all occurrences of $\det$ by $\coldet$)
or, more modestly, to column-commutativity.
When $[a,a] \neq 0$, \reff{eq.proof.lemma.genCayley.1} must be replaced by
\be
   [(\AT B)_{ij}, a_{kl}]  \;=\;  a_{ki} h_{lj}
         \:+\: \sum\limits_m [a_{mi},a_{kl}] \, b_{mj}
   \;,
\ee
and we do not know how to handle the extra terms in
\reff{eq.proof.lemma.genCayley.2}.

4.  For simplicity we have assumed that $s$ is a nonnegative integer,
so that the formulae make sense in arbitrary rings $R$.
But in specific applications we can often allow $s$ to be
an arbitrary real or complex number, or an indeterminate.
For instance, \reff{eq.cor.genCayley.1a}/\reff{eq.cor.genCayley.1b}
make sense in a ring of differential operators
(with $C^\infty$ or analytic coefficients)
over a connected open set $D \subset \R^{n \times n}$ or $\C^{n \times n}$
where $\det X$ is nonvanishing;
here $s$ is an arbitrary real or complex number,
and $(\det X)^s$ denotes any fixed branch
of the corresponding analytic function.
Alternatively, these formulae can be interpreted algebraically,
with $s$ treated as an indeterminate, in a well-known way
\cite[pp.~93--94]{Coutinho_95} \cite[pp.~96~ff.]{Krause_00}.

5.  The Cayley identity \reff{eq.cor.genCayley.2}
has an extension to arbitrary minors,
\be
      \det(\partial_{IJ}) \, (\det X)^s  \;=\;
       s(s+1) \cdots (s+k-1) \, (\det X)^{s-1} \,
       \epsilon(I,J) \, (\det X_{I^c J^c})
 \label{eq.cor.genCayley.2.minors}
\ee
where $|I|=|J|=k$
and $\epsilon(I,J) = (-1)^{\sum_{i \in I} i + \sum_{j \in J} j}$.
Unfortunately we do {\em not}\/ know how to derive
\reff{eq.cor.genCayley.2.minors} from the Capelli identity
(we always seem to get {\em sums}\/ over minors
 rather than one particular minor).
See \cite{CSS_cayley} for alternate combinatorial proofs
of the Cayley identity, which {\em do}\/ include
the all-minors version \reff{eq.cor.genCayley.2.minors}.
\qed

\section*{Acknowledgments}


We wish to thank Luigi Cantini for drawing our attention
to references \cite{Mukhin_06,Chervov_08}.
We also thank an anonymous referee for drawing our attention
to \cite{Chervov_09} and for other helpful comments.

We thank the Isaac Newton Institute for Mathematical Sciences,
University of Cambridge, for support during the programme on
Combinatorics and Statistical Mechanics (January--June 2008),
where this work was (almost) completed.

This research was supported in part by
U.S.\ National Science Foundation grant PHY--0424082.


\begin{thebibliography}{99}


\bibitem{Abdesselam-Crilly-Sokal}  A. Abdesselam, T. Crilly and A.D. Sokal,
   The tangled history of the ``Cayley'' identity
   $\det(\partial) (\det X)^s = s(s+1) \cdots (s+n-1) (\det X)^{s-1}$,
   in preparation.

%

\bibitem{Atiyah_73}  M. Atiyah, R. Bott and V.K. Patodi,
   On the heat equation and the index theorem,
   {\em Invent. Math.}\/ {\bf 19}, 279--330 (1973)
   and erratum {\bf 28}, 277--280 (1975).

\bibitem{Bernstein_72}  I.N. Bern\v{s}te\u{\i}n,
   The analytic continuation of generalized functions with respect to a
   parameter,
   {\em Funkcional. Anal. i Prilo\v{z}en.}\/ {\bf 6}(4), 26--40 (1972)
   [= {\em Funct. Anal. Appl.}\/ {\bf 6}, 273--285 (1972)].


\bibitem{Bjork_79}  J.-E. Bj\"ork, {\em Rings of Differential Operators}\/
    (North-Holland, Amsterdam--Oxford--New York, 1979).

%


\bibitem{Capelli_1882}  A. Capelli, Fondamenti di una teoria generale
   delle forme algebriche, {\em Atti Reale Accad. Lincei, Mem. Classe
   Sci. Fis. Mat. Nat.}\/ (serie 3) {\bf 12}, 529--598 (1882).

\bibitem{Capelli_1887}  A. Capelli, Ueber die Zur\"uckf\"uhrung der
   Cayley'schen Operation $\Omega$ auf gew\"ohnliche Polar-Operationen,
   {\em Math. Annalen}\/ {\bf 29}, 331--338 (1887).

\bibitem{Capelli_1888}  A. Capelli, Ricerca delle operazioni invariantive
   fra pi\`u serie di variabili permutabili con ogni altra operazione
   invariantiva fra le stesse serie,
   {\em Atti Reale Accad. Sci. Fis. Mat. Napoli}\/ (serie 2) {\bf 1},
   1--17 (1888).

\bibitem{Capelli_1890}  A. Capelli, Sur les op\'erations dans
   la th\'eorie des formes alg\'ebriques,
   {\em Math. Annalen}\/ {\bf 37}, 1--37 (1890).

\bibitem{Capelli_02}  A. Capelli, {\em Lezioni sulla Teoria delle Forme
   Algebriche}\/ (Pellerano, Napoli, 1902).


\bibitem{CS_capelli2}  S. Caracciolo and A. Sportiello,
   Noncommutative determinants, Cauchy--Binet formulae,
   and Capelli-type identities.
   II.~Grassmann and quantum oscillator algebra representation,
   in preparation.


\bibitem{CSS_cayley}  S. Caracciolo, A. Sportiello and A.D. Sokal,
   Combinatorial proofs of Cayley-type identities
       for derivatives of determinants and pfaffians, in preparation.

\bibitem{Cartier_69}  P. Cartier and D. Foata,
   {\em Probl\`emes Combinatoires de Commutation et R\'ear\-ran\-ge\-ments}\/,
   Lecture Notes in Mathematics \#85
   (Springer-Verlag, Berlin--Heidel\-berg--New York, 1969).
   Electronic reedition (2006), with three new appendices,
   available on-line at
   \url{http://www.emis.de/journals/SLC/books/cartfoa.html}

\bibitem{Cayley_1846}  A. Cayley, On linear transformations,
   {\em Cambridge and Dublin Math. J.}\/ {\bf 1}, 104--122 (1846).
   [Also in {\em The Collected Mathematical Papers
    of Arthur Cayley}\/ (Cambridge University Press,
    Cambridge, 1889--1897), vol.~1, pp.~95--112.]

\bibitem{Cayley_collected}  A. Cayley, {\em The Collected Mathematical Papers
    of Arthur Cayley}\/, 13 vols. (Cambridge University Press,
    Cambridge, 1889--1897).
    [Also republished by Johnson Reprint Corp., New York, 1963.]

\bibitem{Chervov_08}  A. Chervov and G. Falqui,
   Manin matrices and Talalaev's formula,
   {\em J. Phys. A: Math. Theor.}\/ {\bf 41}, 194006 (2008),
   arXiv:0711.2236 [math.QA] at arXiv.org.

\bibitem{Chervov_09}  A. Chervov, G. Falqui and V. Rubtsov,
   Algebraic properties of Manin matrices 1,
   arXiv:0901.0235 [math.QA] at arXiv.org.

%

\bibitem{Coutinho_95}  S.C. Coutinho,
   {\em A Primer of Algebraic $D$-Modules}\/,
   London Mathematical Society Student Texts \#33
   (Cambridge University Press, Cambridge, 1995).


\bibitem{Dolgachev_03}  I. Dolgachev, {\em Lectures on Invariant Theory}\/,
   London Mathematical Society Lecture Note Series \#296
   (Cambridge University Press, Cambridge, 2003).



\bibitem{Foata_79}  D. Foata, A noncommutative version of the matrix
   inversion formula,
   {\em Adv. Math.}\/ {\bf 31}, 330--349 (1979). 

\bibitem{Foata_94}  D. Foata and D. Zeilberger, Combinatorial proofs of
   Capelli's and Turnbull's identities from classical invariant theory,
   {\em Electron. J. Combin.} {\bf 1}, \#R1 (1994).

\bibitem{Fulton_91}  W. Fulton and J. Harris,
   {\em Representation Theory: A First Course}\/
   (Springer-Verlag, New York--Heidelberg, 1991), Appendix F.



%
%
%
%

\bibitem{Howe_89}  R. Howe, Remarks on classical invariant theory,
   {\em Trans. Amer. Math. Soc.}\/ {\bf 313}, 539--570 (1989)
   and erratum {\bf 318}, 823 (1990).

\bibitem{Howe_91}  R. Howe and T. Umeda, The Capelli identity,
   the double commutant theorem, and multiplicity-free actions,
   {\em Math. Ann.}\/ {\bf 290}, 565--619 (1991).

\bibitem{Itoh_00}  M. Itoh, Capelli elements for the orthogonal Lie algebras,
   {\em J. Lie Theory}\/ {\bf 10}, 463--489 (2000).

\bibitem{Itoh_04}  M. Itoh, Capelli identities for the dual pair
   $(O_M, Sp_N)$, {\em Math. Z.}\/ {\bf 246}, 125--154 (2004).

\bibitem{Itoh_05}  M. Itoh, Capelli identities for reductive dual pairs,
    {\em Adv. Math.}\/ {\bf 194}, 345--397 (2005).

\bibitem{Itoh_07}  M. Itoh, Two determinants in the universal enveloping
    algebras of the orthogonal Lie algebras,
    {\em J. Algebra}\/ {\bf 314}, 479--506 (2007).

\bibitem{Itoh-Umeda_01}  M. Itoh and T. Umeda, On central elements in the
    universal enveloping algebras of the orthogonal Lie algebras,
    {\em Compositio Math.}\/ {\bf 127}, 333--359 (2001).

\bibitem{Kinoshita_02}  K. Kinoshita and M. Wakayama,
   Explicit Capelli identities for skew symmetric matrices,
   {\em Proc. Edinburgh Math. Soc.}\/ {\bf 45}, 449--465 (2002).


\bibitem{Konvalinka_08a} M. Konvalinka, Non-commutative Sylvester's
   determinantal identity,
   {\em Electron. J. Combin.}\/ {\bf 14}, \#R42 (2008),
   arxiv:math/0703213 [math.CO] at arXiv.org.

\bibitem{Konvalinka_08b}  M. Konvalinka, An inverse matrix formula in the
   right-quantum algebra,
   {\em Electron. J. Combin.}\/ {\bf 15}, \#R23 (2008),
   arxiv:math/0703203 [math.CO] at arXiv.org.

\bibitem{Konvalinka_07}  M. Konvalinka and I. Pak,
   Non-commutative extensions of the MacMahon Master Theorem,
   {\em Adv. Math.}\/ {\bf 216}, 29--61 (2007),
   arXiv:math/0607737 [math.CO] at arXiv.org.

\bibitem{Kostant_91}  B. Kostant and S. Sahi, The Capelli identity,
   tube domains, and the generalized Laplace transform,
   {\em Adv. Math.}\/ {\bf 87}, 71--92 (1991).

\bibitem{Kostant_93}  B. Kostant and S. Sahi, Jordan algebras and
   Capelli identities, {\em Invent. Math.}\/ {\bf 112}, 657--664 (1993).

\bibitem{Kraft_96}  H. Kraft and C. Procesi, {\em Classical Invariant Theory:
   A Primer}\/ (preliminary version, July 1996),
   available on-line at
   \url{http://www.math.unibas.ch/~kraft/Papers/KP-Primer.pdf}

\bibitem{Krause_00}  G.R. Krause and T.H. Lenagan,
   {\em Growth of Algebras and Gelfand--Kirillov Dimension}\/,
   rev.\ ed.\ 
   (American Mathematical Society, Providence RI, 2000).

\bibitem{Lalonde_96}  P. Lalonde, A non-commutative version of
    Jacobi's equality on the cofactors of a matrix,
    {\em Discrete Math.}\/ {\bf 158}, 161--172 (1996).

%


\bibitem{Manin_87}  Yu. I. Manin, Some remarks on Koszul algebras and
   quantum groups, {\em Ann. Inst. Fourier (Grenoble)}\/ {\bf 37},
   191--205 (1987).

\bibitem{Manin_88}  Yu.I. Manin, {\em Quantum Groups and Non-Commutative
   Geometry}\/ $\;$
   (Centre de Recherches Math\'ematiques, Universit\'e de Montr\'eal, 1988).
   
\bibitem{Manin_91}  Y.I. Manin, {\em Topics in Noncommutative Geometry}\/
   (Princeton University Press, Princeton, NJ, 1991).


\bibitem{Molev_07}  A. Molev, {\em Yangians and Classical Lie Algebras}\/,
   Mathematical Surveys and Monographs \#143
   (American Mathematical Society, Providence, RI, 2007).

\bibitem{Molev_99}  A. Molev and M. Nazarov, Capelli identities for
    classical Lie algebras, {\em Math. Ann.}\/ {\bf 313}, 315--357 (1999).

%

\bibitem{Mukhin_06}  E. Mukhin, V. Tarasov and A. Varchenko,
   A generalization of the Capelli identity,
   preprint (October 2006),
   arXiv:math/0610799 [math.QA] at arXiv.org.

%
%

\bibitem{Nazarov_97}  M. Nazarov, Capelli identities for Lie superalgebras,
   {\em Ann. Sci. \'Ecole Norm. Sup.}\/ {\bf 30}, 847--872 (1997).

\bibitem{Okounkov_96a}  A. Okounkov, Quantum immanants and higher Capelli
   identities, {\em Transform. Groups}\/ {\bf 1}, 99--126 (1996).

\bibitem{Okounkov_96b}  A. Okounkov, Young basis, Wick formula,
   and higher Capelli identities, {\em Internat. Math. Res. Notices}\/
   {\bf 1996}, no. 17, 817--839.

\bibitem{Olver_99}  P.J. Olver, {\em Classical Invariant Theory}\/,
   London Mathematical Society Student Texts \#44
   (Cambridge University Press, Cambridge, 1999).

\bibitem{Nazarov_98}  M. Nazarov, Yangians and Capelli identities
   (Kirillov's seminar on representation theory),
   {\em Amer. Math. Soc. Transl. Ser. 2}\/ {\bf 181}, 139--163 (1998).

\bibitem{Noumi_94}  M. Noumi, T. Umeda and M. Wakayama,
   A quantum analogue of the Capelli identity and an elementary differential
   calculus on ${\rm GL}\sb q(n)$,
   {\em Duke Math. J.}\/ {\bf 76}, 567--594 (1994).

\bibitem{Noumi_96}  M. Noumi, T. Umeda and M. Wakayama,
   Dual pairs, spherical harmonics and a Capelli identity in quantum group
   theory, {\em Compositio Math.}\/ {\bf 104}, 227--277 (1996).



%

\bibitem{Schur_68}  I. Schur, {\em Vorlesungen \"uber Invariantentheorie}\/
   (Springer-Verlag, Berlin--Heidelberg--New York, 1968).


%


\bibitem{Stein_82} P.R. Stein, On an identity from classical invariant theory,
   {\em Lin. Multilin. Alg.}\/ {\bf 11}, 39--44 (1982).

%
%

\bibitem{Turnbull_48}  H.W. Turnbull, Symmetric determinants and the
   Cayley and Capelli operators, {\em Proc. Edinburgh Math. Soc.}\/
   (2){\bf 8}, 76--86 (1948).


\bibitem{Umeda_98}  T. Umeda, The Capelli identities, a century after,
   {\em Amer. Math. Soc. Transl. Ser. 2}\/ {\bf 183}, 51--78 (1998).

\bibitem{Umeda_08}  T. Umeda, On the proof of the Capelli identities,
   {\em Funkcial. Ekvac.}\/ {\bf 51}, 1--15 (2008).


\bibitem{Wachi_06}  A. Wachi, Central elements in the universal enveloping
    algebras for the split realization of the orthogonal Lie algebras,
    {\em Lett. Math. Phys.}\/ {\bf 77}, 155--168 (2006).

\bibitem{Wallace_53}  A.H. Wallace, A note on the Capelli operators
   associated with a symmetric matrix, {\em Proc. Edinburgh Math. Soc.}\/
   (2){\bf 9}, 7--12 (1953).


\bibitem{Weyl_46}  H. Weyl, {\em The Classical Groups, Their Invariants
    and Representations}\/, 2nd ed.\ 
    (Princeton University Press, Princeton NJ, 1946).




\end{thebibliography}
\end{document}